\documentclass[11pt,a4paper]{article}
\pdfoutput=1
\RequirePackage{amsmath}%
\usepackage{amsfonts}%
\usepackage{amssymb}%
\usepackage{graphicx}%
\usepackage{url}%
\usepackage{color}

\newtheorem{theorem}{Theorem}

\newtheorem{lem}{Lemma}

\newtheorem{remark}{Remark}

\usepackage[autolinebreaks]{mcode}
\usepackage{url}
\usepackage[title]{appendix}%

\usepackage{floatflt,graphicx,graphics,subfig}

\DeclareMathOperator{\capa}{\mathrm{cap}}
\let\Re\relax
\DeclareMathOperator{\Re}{\mathrm{Re}}
\let\Im\relax
\DeclareMathOperator{\Im}{\mathrm{Im}}

\renewcommand{\i}{\mathrm{i}}
\newcommand{\bI}{{\bf I}}
\newcommand{\bM}{{\bf M}}
\newcommand{\bN}{{\bf N}}

\newcommand{\M}{\mathsf{M}}
\newcommand{\K}{\mathcal{K}}

\newcommand{\cT}{\mathcal{T}}
\newcommand{\B}{\mathbb{B}^2}
\DeclareMathOperator{\sh}{\mathrm{sh}}
\renewcommand{\th}{\,\textnormal{th}}
\newcommand{\BB}{\mathbb{B}^2}
\newcommand{\arth}{\,\textnormal{arth}}

\newcommand{\eE}{\varepsilon}
\newcommand{\osc}{\mathrm{osc}}

\begin{document}
\title{Constrained maximization of conformal capacity}
	\author{Harri Hakula$^{\rm a}$, Mohamed M.S. Nasser$^{\rm b}$, Matti Vuorinen$^{\rm c}$}
	
	\date{}
	\maketitle
 	
	\vskip-0.8cm %
	\centerline{$^{\rm a}$Department of Mathematics and System Analysis, P.O. Box 11100,} %
\centerline{FI--00076 Aalto, Finland}%
	\centerline{\tt harri.hakula@aalto.fi}%
	\centerline{$^{\rm b}$Department of Mathematics, Statistics \& Physics, Wichita State University,} %
	\centerline{Wichita, KS 67260-0033, USA}%
	\centerline{\tt mms.nasser@wichita.edu}%
	\centerline{$^{\rm c}$Department of Mathematics and Statistics, University of Turku, Turku, Finland} %
	\centerline{\tt vuorinen@utu.fi}%

\begin{abstract}
We consider constellations of disks which are unions of disjoint hyperbolic disks in the unit disk with fixed radii and unfixed centers. 
We study the problem of maximizing the conformal capacity of a constellation
with a fixed number of disks under constraints on the centers in two cases.
In the first case the constraint is that the centers are at most at distance $R \in(0,1)$ from the origin and in the second case it is required that the centers are on the subsegment $[-R,R]$ of a diameter of the unit disk. 
We study also similar types of constellations with hyperbolic segments instead of the hyperbolic disks.
Our computational experiments suggest that a dispersion phenomenon occurs: the disks/segments go as close to the unit circle as possible under these constraints and stay as far as possible from each other. 
The computation of capacity reduces to the Dirichlet problem for the Laplace equation which we solve using two methods: a fast boundary integral equation method and a high-order finite element method.
\end{abstract}

\begin{center}
\begin{quotation}
{\noindent {{\bf Keywords}.\;\; Conformal capacity, condenser, hyperbolic geometry, boundary integral equation method, $hp$-FEM, Dirichlet problem}%
}%
\end{quotation}
\end{center}

\begin{center}
\begin{quotation}
{\noindent {{\bf MSC}.\;\; 65E05, 31A15, 30C85 }}%
\end{quotation}
\end{center}

\section{Introduction}

In physics and chemistry there exist vast number of problems that involve interactions between
multiple bodies or particles, such as the $n$-body problem of celestial mechanics and the
many-body problem in quantum physics. Such interactions can be pairwise, such as forces, or
complicated effects through some aggregate fields. Typically one is interested in configurations
that imply some extremal state, for instance, minimum or maximum energy ones. These
extremal configurations often have geometric features, such as symmetries, and give rise to many
packing problems.  One example of such problems is the modeling of intermolecular forces in simple systems
using Lennard-Jones potential, which incorporates both
short-range repulsion and long-range attraction \cite{do,wa}.
The transition between repulsive and attractive forces
occurs at a specific equilibrium distance defining the most
stable configuration where the atoms forming the molecule prefer to remain.

Here we study
condensers of the form $(\mathbb{B}^2, E)$ where for a fixed integer $m\ge 2,$ the
set $E$ is a union of 
finitely many disjoint closed disks $E_1,\ldots,E_m$ in the unit  disk $\mathbb{B}^2$  with fixed hyperbolic
radii. We call such a collection of sets, 
or for that matter also the set $E,$ {\it a constellation} of disks.
Note that these disks are allowed to move: only the hyperbolic radii are fixed, but
the centers are not. The location of the disks within a constellation has
a strong influence on the numerical value of the conformal capacity.
Our goal is to study extremal problems for the conformal capacity of condensers 
of the form $(\mathbb{B}^2,E)$ where  $\mathbb{B}^2$ is the unit disk and  
the set $E \subset \mathbb{B}^2$ is a  constellation of hyperbolic disks.

Classical results show that
applying a geometric transformation, so called symmetrization, on a compact set $F \subset \mathbb{B}^2,$
the new symmetrized set $F^*$ exhibits some symmetry and what is relevant here,
the new set $F^*$ provides a lower bound for the conformal capacity~\cite{ba,du,g,ps}
\begin{equation} \label{symbd}
	{\rm cap} (\mathbb{B}^2,F) \ge {\rm cap} (\mathbb{B}^2,F^*).
\end{equation}
Equality holds here if $F= F^*.$ Due to the conformal 
invariance of the conformal capacity,
in the case $F \neq F^*$, the lower bound \eqref{symbd} can be
improved at least in the case when $F$ is a constellation of separate 
hyperbolic disks. In the recent paper \cite{hnv}, the conformally invariant
hyperbolic geometry was used as a key tool to refine \eqref{symbd}.

We study here a reverse problem, maximization of conformal capacity.
Symmetrization methods applied to a set $E$ often reduce the distances
between the points in $E$ while some set functional like the area
remains invariant. In the maximization process a reverse phenomenon can
be naturally expected. Some results have  been reported in
literature~\cite{bsv,knv}. Since the publication of the classical
monograph~\cite{ps}, many authors have studied extremal problems from
the point of view of potential theory~\cite{ba,bhs,dfhsz,du,sz}.

It seems natural to study the problem of capacity maximization of a
constellation under suitable constraints. 
We study two cases: (i)~the
centers of the disks of the constellation are contained in a subdisk, $E
\subset B^2(0,R)$, (ii)~the centers of the disks of the constellation
are on a symmetric subsegment $[-R,R]$ on a the diameter $(-1,1)$ of 
$\mathbb{B}^2$. In both cases our simulations suggest that some kind of
a maximal dispersion phenomenon occurs similar to what is observed
for instance in chemistry: the disks increase their mutual
``social distances'' and, at the same time, have a tendency to move as
close to the unit circle as the constraints permit. This phenomenon is
the reason why a constraint for the centers of the disks is natural:
without such a constraint, during the maximization process, the disks
could go arbitrarily close to the unit circle and become ``invisible'',
their Euclidean diameters would become arbitrarily small. Similar
results are obtained when the hyperbolic disks in the above
constellations are replaced by hyperbolic segments of fixed lengths
such as radial hyperbolic segments
and hyperbolic segments on the real line. 
Based on the numerical experiments we can predict the maximal capacities
of all constellations exhibiting symmetry, that is, when the
constellation in its extremal configuration can be divided into identical
compartments of the unit disk.

\subsection{Illustrative Example}
\begin{figure}
	\centering
	\includegraphics[width=0.475\textwidth]{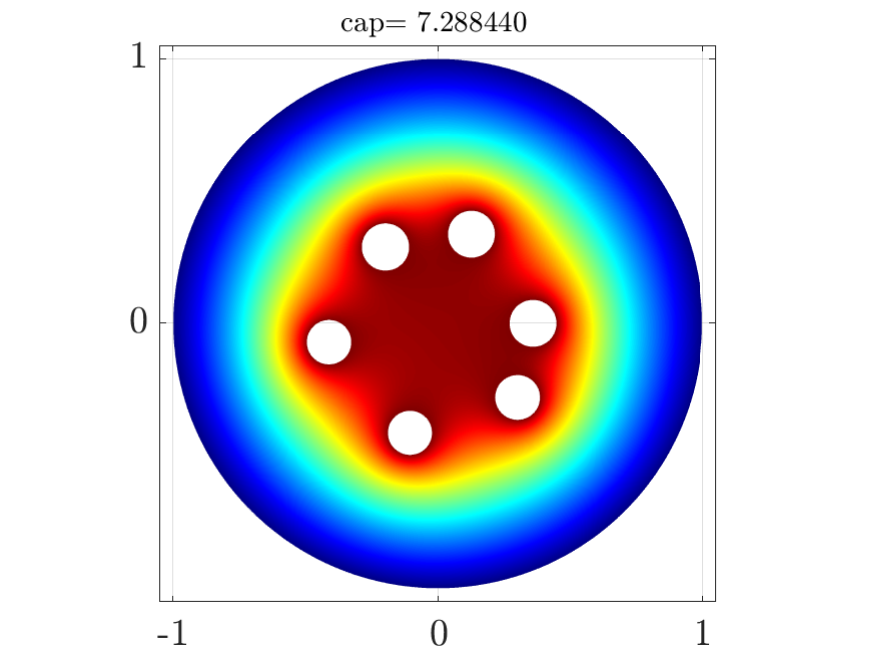} \hfill
	\includegraphics[width=0.475\textwidth]{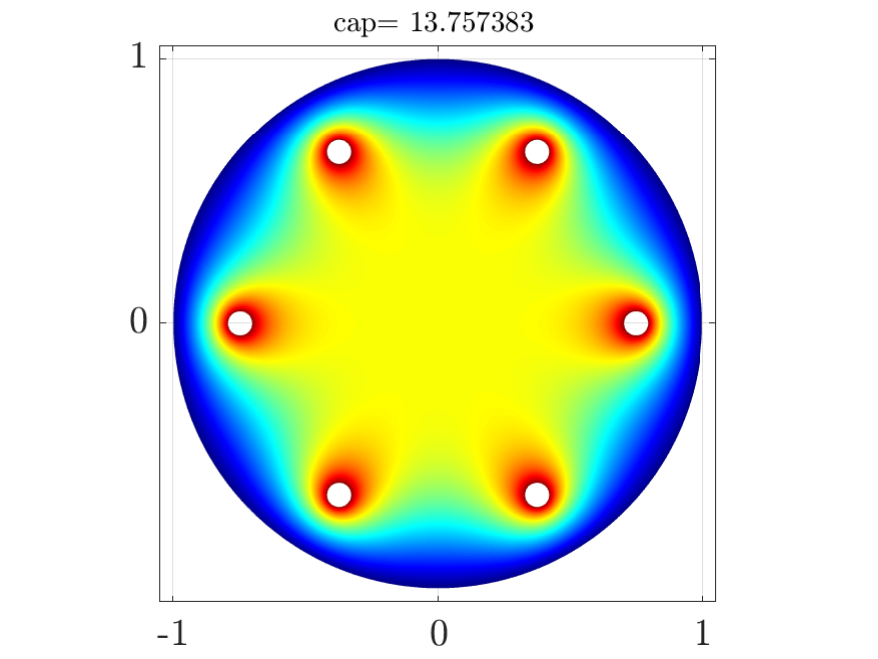} \\
	\includegraphics[width=0.475\textwidth]{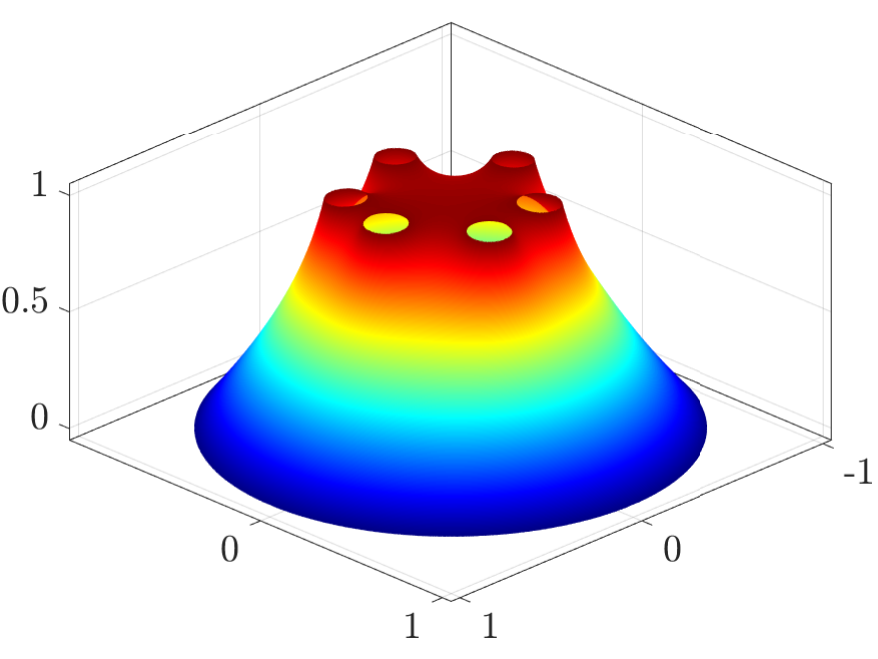} \hfill
	\includegraphics[width=0.475\textwidth]{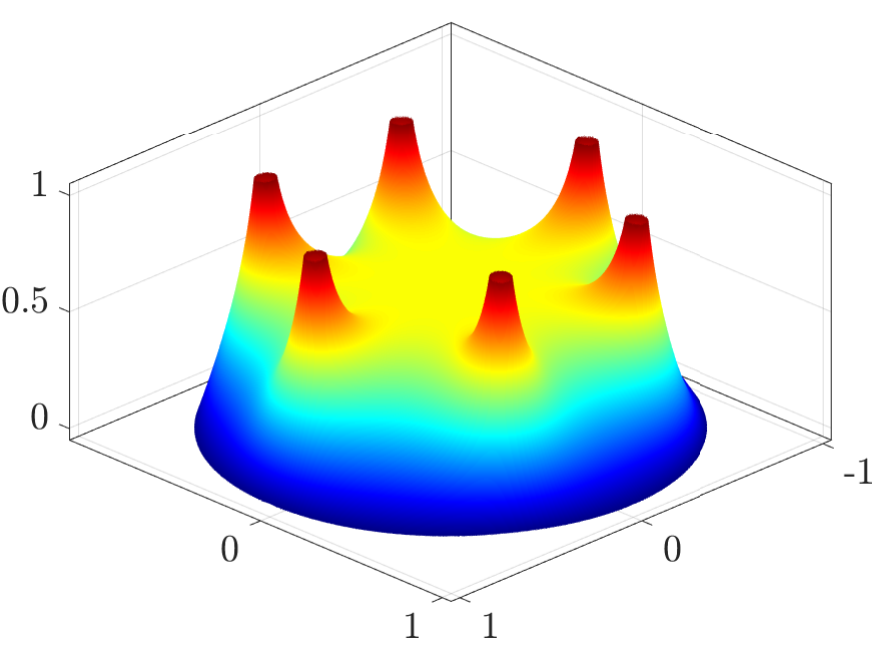} 
	\caption{Surface plots of the potentials for the six hyperbolic disks with equal hyperbolic radii $0.2$. The centers of the disks are at the initial positions (left) and at the positions that maximize the capacity (right).}\label{fig:6Dint}
\end{figure}
Consider a constellation $E$ of six disks with equal hyperbolic radii 
whose centers are contained within a subdisk, 
$E \subset B^2(0,R)$. The task is to find a configuration with maximal capacity.
This is illustrated in Figure~\ref{fig:6Dint}. Maximizing the capacity is equivalent to
 maximizing the $L^2$-norm of the gradient of the potential $u$ which is a solution of
the Dirichlet problem for the Laplace equation
where the Dirichlet boundary values at the points of the unit circle are zero 
and at the boundary points of each disk equal to one.
(For the formal definition of the capacity, see Section~\ref{sec:capacity} below.)

Analyzing the results of optimization one can observe that, regardless of the initial configuration, two types of features of the final configuration of the constellation.
First, the  disks move towards the unit circle, as close as the constraint permits, and second, their mutual hyperbolic distances appear to become maximal,
resulting in a symmetric configuration. This final configuration demonstrates what
is referred to above as a maximal dispersion phenomenon.

\subsection{Organization}
The paper is organized as follows: Section 2 contains
preliminary information about hyperbolic geometry, conformal capacity, and special
functions to be used in the later sections. Section 3 is a description of our two
computational methods, the $hp-$FEM and the boundary integral equation method.
Our experimental discoveries are confirmed by these two methods. Section 4 presents our computational work on the disk constellations. 
Section 5 presents similar results, but now in place of the hyperbolic disks we have hyperbolic segments with fixed lengths.
Section~6 draws the conclusions of our work and suggests problems for new research. 

\section{Preliminaries}
In this section we recall some facts from hyperbolic geometry and special functions related to conformal capacity of canonical condensers.

\subsection{Hyperbolic geometry}\label{hg}
We recall some basic formulas and notation for hyperbolic geometry from~\cite{b}.
The Euclidean balls with center $x \in {\mathbb{R}^2}$ and radius $r>0$ are denoted $B^2(x,r)$ and its boundary sphere is $S(x,r)$. For brevity we write $ {\mathbb{B}^2}= B^2(0,1)$.
For $a,b\in \mathbb{B}^2$, the hyperbolic distance $\rho(a,b)$ between $a$ and $b$ is defined via the formula
\begin{equation}\label{rhoB}
{\rm sh} \frac{\rho(a,b)}{2}=
\frac{|a-b|}{\sqrt{(1-|a|^2)(1-|b|^2)}}.
\end{equation}
The hyperbolic disk with center $x \in {\mathbb{B}^2} $ and radius $R>0$
is $B_{\rho}(x,R)=\{ z: \rho(x,z)<R \}.$ We often use the
connection between the hyperbolic disk and Euclidean disk
\begin{equation}\label{hkv4.20}
\begin{cases}
B_\rho(x,R)=B^2(y,r)\;,&\\
\noalign{\vskip5pt}
{\displaystyle y={x(1-t^2) \over 1-|x|^2t^2}\;,\;\;
r={(1-|x|^2)t \over 1-|x|^2t^2}\;,\;\;t={\rm th}(R/2)\;,}&
\end{cases}
\end{equation}
We use the notation ${\rm sh}$ and ${\rm th}$ for the hyperbolic sine and the hyperbolic tangent, respectively.  Their inverse functions are ${\rm arsh}$ and ${\rm arth}$.

\subsection{Special functions}
	For $|z|<1$, 
	the Gaussian hypergeometric function is defined by the equality
	$$
	{}_2F_1(a,b;c;z)=\sum_{n=1}^\infty
	\frac{(a)_n(b)_n}{(c)_n}\,\frac{z^n}{n!},
	$$
	where $(q)_n$ denotes the Pochhammer symbol, i.e.
	$(q)_n=q(q+1)\ldots(q+(n-1))$ for every natural $n$ and $(q)_0=1$~\cite{as}.
	
	The complete elliptic integral of the first kind
	\begin{equation}\label{ek}
		{\K}
		(r)=\int_0^1\frac{dt}{\sqrt{(1-t^2)(1-r^2t^2)}}\,,\quad r \in (0,1),
	\end{equation}
	is, in fact, a special case of the Gaussian hypergeometric
	function; we have
	$$
	{\K}
	(r)=\frac{\pi}{2}\,\,{}_2F_1({\textstyle\frac{1}{2}},{\textstyle\frac{1}{2}};1;r^2).
	$$
	
	\noindent 
	The decreasing homeomorphism $\mu:(0,1)\to (0,\infty)$
	$$
	\mu(r)=\frac{\pi}{2} \frac{{\K}(\sqrt{1-r^2})}{{\K}(r)}\,,\quad
	0<r<1,
	$$
	is recurrent in the study of conformal invariants.

\subsection{Condenser capacity}\label{sec:capacity}
A \emph{condenser} is a pair $(G,E)$, where $G \subset \mathbb{B}^2$ 
is a domain and $E$ is a compact non-empty subset of $G$. 
The \emph{conformal capacity} of this condenser is defined as  
\cite{du,GoRe,hkv}
\begin{align}\label{def_condensercap}
{\rm cap}(G,E)=\inf_{u\in A}\int_{G}|\nabla u|^2 dm,
\end{align}
where $A$ is the class of $C^\infty_0(G)$ functions $u: G\to[0,\infty)$ 
with  $u(x) \ge 1$ for all $x \in E$ and $dm$ is the $2$-dimensional 
Lebesgue measure. 
Here we assume that $G=\mathbb{B}^2$ is the unit disk and $E=\cup_{j=1}^{m}E_j$ where $E_1,\ldots,E_m$ are $m$ compact disjoint non-empty subsets of the unit disk $\mathbb{B}^2$ such that $\partial E_1, \ldots, \partial E_m$ are smooth Jordan curves. 
Hence $\Omega=G\backslash E$ is a multiply connected domain of connectivity $m+1$ and the infimum in~\eqref{def_condensercap} is attained by a harmonic function $u$. This extremal function $u$ is the unique solution of the Laplace equation in $\Omega$ with boundary values equal $u=1$ on $E$ and $u=0$ on $\partial G$~\cite{du}.
The capacity can be then expressed in terms of the extremal function $u$ as
\begin{equation}\label{eq:cap}
	{\rm cap}(G,E)=\iint \limits_{\Omega}|\nabla u|^2 dxdy,
\end{equation}
which, using Green's formula~\cite[p. 4]{du}, implies that
\begin{equation}\label{eq:cap2}
	{\rm cap}(G,E)=\int_{\partial\Omega} u\frac{\partial u}{\partial {\bf n}}ds
\end{equation}
where $\partial u/\partial {\bf n}$ denotes the directional derivative of $u$ along the outward normal. Since $u=0$ on $\partial G$ and $u=1$ on $\partial E_k$, we have  
\begin{equation}\label{eq:cap3}
	{\rm cap}(G,E)=\sum_{k=1}^m b_k
\end{equation}
where 
\begin{equation}\label{eq:cap-bk}
	b_k=\int_{\partial E_k} \frac{\partial u}{\partial {\bf n}}ds, \quad k=1,2,\ldots,m.
\end{equation}
Thus, the constant $b_k$ can be considered as the contribution of the compact set $E_k$ to the capacity ${\rm cap}(G,E)$,  for $k=1,2,\ldots,m$.
Since the Dirichlet integral is conformally invariant, the cases for which $\partial E_1, \ldots, \partial E_m$ are rectilinear slits can be handled with the help of auxiliary conformal mappings which transform the slits to smooth curves.

The conformal capacity of a condenser is one 
of the key notions of potential theory of elliptic partial 
differential equations \cite{GoRe,hkm} and it has numerous 
applications to geometric function theory, both in the plane and 
in higher dimensions, \cite{du,  GoRe, hkv, hkm}.

	Numerous variants of the definition \eqref{def_condensercap} of 
	capacity are given in \cite{GoRe,hkv}. 
	For instance
	\begin{align} \label{ziem}
		{\rm cap}(G,E)=\M(\Delta(E,\partial G;G)),    
	\end{align}
	where $\Delta(E,\partial G;G)$ is the family of all curves joining $E$ 
	with the boundary $\partial G$ in the domain $G$ and $\M$ stands 
	for the modulus of a curve family \cite[Ch~7]{hkv}.
	A fundamental fact is {\it subadditivity}: if $E= \cup_{j=1}^m E_j$ where
	$E_j \subset \mathbb{B}^2$ for all $j,$ then
	\begin{equation} \label{subadd}
		{\rm cap}(\mathbb{B}^2,E)\le \sum_{j=1}^m{\rm cap}(\mathbb{B}^2,E_j).
	\end{equation}
	For the basic facts about capacities and moduli, the reader is 
	referred to \cite{du, GoRe, hkv, hkm}.

	The exact value of the capacity
	is known only in a handful of special cases.  For instance, the
	capacity $\gamma_2(r)$ of the Gr\"otzsch condenser $(\mathbb{B}^2,[0,r])$ 
	can be expressed as
	$$
	\gamma_2(r)=2\pi/\mu(r).
	$$
	The capacity of an annulus is also known by the next lemma and~\eqref{ziem}.

\begin{lem}\label{cgqm_5.14}\cite[(7.3), p. 107]{hkv}\newline
	(1) If $0<a<b$ and $D=\overline{B}^2(0,b)\backslash B^2(0,a)$,
	\begin{align*}
		\M(\Delta(S(0,a),S(0,b);D))=2\pi/\log({b}/{a}).
	\end{align*}
	(2) If $R>0$ then for $x \in \mathbb{B}^2$ 
	\begin{align*}
		\M(\Delta(S(0,1),B_{\rho}(x,R);\mathbb{B}^2))=
		2\pi/\log({1}/{{\rm th} (R/2)}) \,.
	\end{align*}
\end{lem}

\subsection{Modulus of a quadrilateral.}\label{quadri}
A \textit{ quadrilateral} is a Jordan domain  $Q$ on the
complex plane with four marked points $z_1$, $z_2$, $z_3$, $z_4$ on
the boundary. We assume that the points define positive orientation.   
By Riemann's mapping theorem, there is a conformal mapping of $Q$ onto a
rectangle  with vertices $1,1+hi,hi,0$, $h>0$,  such that the vertices
of $Q$ correspond to the vertices of the rectangle. Then the value $h$
 is called the conformal modulus of $Q$:
$$
h ={\rm mod} Q \equiv \M(\Delta([0,1],[ih, 1+ih]; Q)) .
$$
An alternative method to find the modulus  is to solve  the 
following Dirichlet-Neumann boundary value problem for the Laplace equation \cite{du}. 
Suppose that $\partial Q = \cup_{k=1}^4\partial Q_k$; all the four boundary arcs between vertices are
assumed to be non-degenerate. This problem is
\begin{equation*} 
\left\{\begin{matrix}
    \Delta u& =&\ 0,& \text{on}\ &{\ Q,} \\
    u& =&\ 1,&  \text{on}\ &{\partial Q_1 = (z_1,z_2),}\\
    u& =&\ 0, & \text{on}\ &{\partial Q_3 = (z_3,z_4),}\\
    \partial u/\partial n&  =&\ 0,& \text{on}\ &{\partial Q_2 = (z_2,z_3),}\\
    \partial u/\partial n&  =&\ 0,& \text{on}\ &{\partial Q_4 = (z_4,z_1).}\\
\end{matrix}\right.
\end{equation*}
In terms of a solution function $u$ to the above problem,  the modulus can be
computed  as \[
\int\!\!\int_Q |\nabla u|^2 dm \,.
\]

\subsection{Symmetric compartments.} 
\begin{figure}
	\centering
	\subfloat{\includegraphics[width=0.5\linewidth]{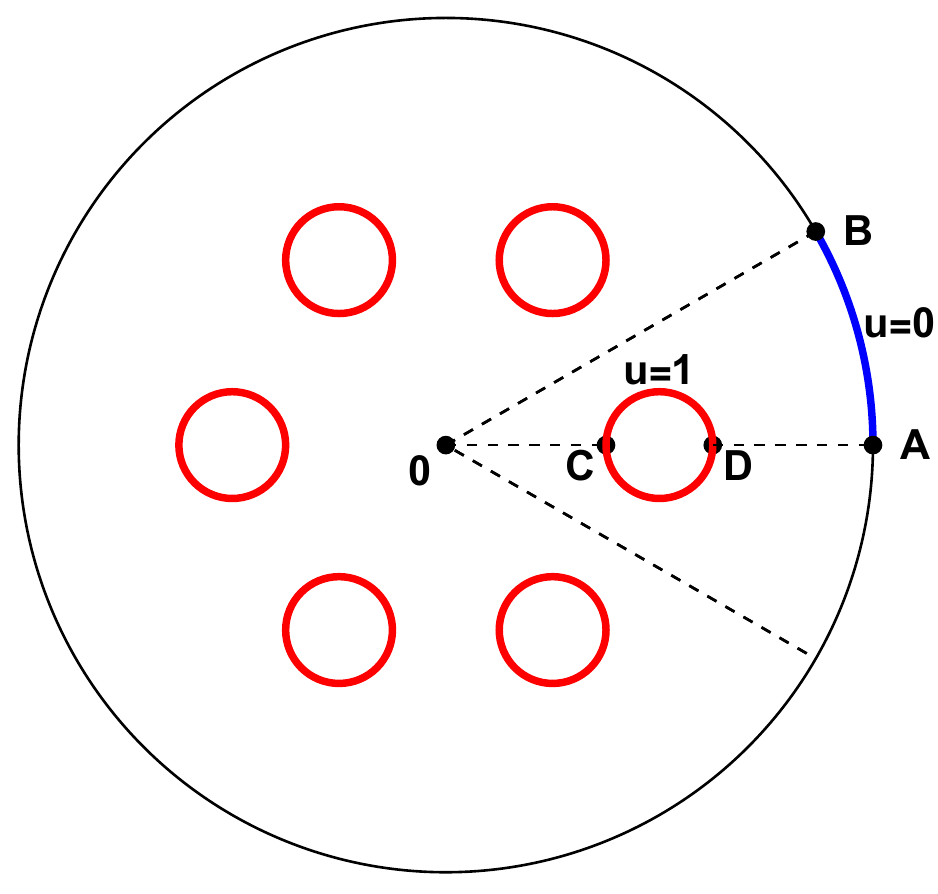}}
	\caption{Sectorial symmetry. The capacity of each sector or compartment can be computed
	separately. In the case of a constellation with identical elements, each compartment can be divided into
	two quadrilaterals symmetric with respect to the radius bisecting the sector.}
	\label{fig:figCompart}
\end{figure}

If the constellation has suitable symmetry, like in the case of disks with equal size in Figure \ref{fig:6Dint}, one can subdivide
the domain into non-overlapping compartments, carry out the computation in each
compartment. Due to symmetry, the result is equal for each compartment and we may thus
reduce the computational load by a factor which in the case of $m$ sets is roughly 
$1/(2m).$ To explain this idea further, we use moduli of quadrilaterals and refer
to Figure \ref{fig:figCompart}.

The disks $E_j, j=1,...,m,$ of the constellation have equal radii and
their centers are $r\exp(\i 2k \pi/m), k=0,1,...,m, m=6,$ for some $r
\in(0, 1)$ such that these $m$ disks and the unit circle are not
overlapping. 
The first compartment  is $\{ - \pi/m< \arg(z)<\pi/m\}$ and the first
disk of the constellation is symmetric with respect to the real axis. In the
upper half of the first compartment we consider the quadrilateral $Q$
with vertices $A,B,C,D$ and the associated
Dirichlet-Neumann boundary value problem with $u(z) =1$ on the red semicircle and $u= 0$
on the blue subarc of the unit circle and $\partial u/\partial n = 0 $ on the dashed sides of the upper half
of the compartment. Our conclusion is that
\begin{equation} \label{capCompart}
{\rm cap}(\mathbb{B}^2, E) = 2 m  \, {\rm mod} Q \,,
\end{equation}
where $E=\cup_{j=1}^m E_j$.
We can use this relation to predict the limiting value of the maximization of the capacity, e.g., of the constellation in Figure~\ref{fig:6Dint}. 
We present in Table~\ref{tab:Q} the numerical values of ${\rm cap}(\mathbb{B}^2, E)$ obtained using the BIE method and numerical values of $2 m  \, {\rm mod} Q$ obtained using the FEM for several values of $m$. We assume here that the disks $E_1,\ldots,E_m$ are Euclidean disks such that the center and radius of the disk $E_j$ is $0.5e^{\i 2(j-1)\pi/m}$ and $0.1$, $j=1,\ldots,m$.
The vertices of the quadrilateral in Figure~\ref{fig:figCompart} are then $A=1$, $B=e^{\i \pi/m}$, $C=0.4$, and $D=0.6$.

\begin{table}
	\caption{The numerical values of ${\rm cap}(\mathbb{B}^2, E)$ and $2 m  \, {\rm mod} Q$.}
	\label{tab:Q}
	\centering
	\begin{tabular}{lccc} 
		\hline
		$m$ & ${\rm cap}(\mathbb{B}^2, E)$ & $2 m  \, {\rm mod} Q$ & $|{\rm cap}(\mathbb{B}^2,E)-2m\,{\rm mod}Q|$ \\ \hline
		$5$ & $9.47487674904924$           &    $9.47487674904923$  & $1\times 10^{-14}$   \\
		$6$ & $10.0486182568334$           &    $10.04861825683339$ & $1\times 10^{-14}$   \\
		$7$ & $10.4636668610180$           &    $10.46366686101804$ & $4\times 10^{-14}$   \\
		$8$ & $10.7735173309461$           &    $10.77351733094614$ & $4\times 10^{-14}$   \\
		\hline
	\end{tabular}
\end{table}

\section{Methods}

In this section the numerical methods used in the numerical experiments
are briefly described. The capacities of constellations are computed using 
the boundary integral equation with
the generalized Neumann kernel method (BIE) impleneted in MATLAB
and 
the $hp$-version of the finite element method (FEM) implemented in Mathematica.
The maximization problems are computed using the interior-point method
as implemented in MATLAB and Mathematica.

In any numerical study the questions of validation and verification need to be addressed.
The Dirichlet problem~\eqref{def_condensercap}
is one of the primary numerical model problems, 
therefore any standard solution technique can be viewed as having been validated.
For example, for the constellations of disks case, i.e., when the domain $\Omega$ is a multiply connected circular domain of connectivity $m+1$, the Dirichlet problem~\eqref{def_condensercap} can be solved using linear combinations of $m$ multi-valued analytic functions $v_1(z),\ldots,v_m(z)$ relevant to a special transcendental function known as the Schottky-Klein prime function associated with the circular domain $\Omega$~\cite{CrB,ckgn,MC23}. Two numerical methods have been presented in~\cite{ckgn} for computing the $m$ functions $v_1(z),\ldots,v_m(z)$. One of these methods is based on using the boundary integral equation with the generalized Neumann kernel that will be used in this paper.
Computing these $m$ functions $v_1(z),\ldots,v_m(z)$ requires solving $m$ integral equations. Thus, in principle, the computational cost of computing these $m$ functions using the integral equation and then using these $m$ functions to compute the capacity is the same as the computational cost of using the integral equation to compute the capacity using the method discussed in this paper.
Furthermore, explicit formulas for estimating the conformal capacity were derived in~\cite{MC23} by using the matching approach. 

In the class of problems considered in this paper, verification follows through using two numerical methods. Namely, a high-order finite element method ($hp$-FEM) and a boundary integral equation method (BIE). This is illustrated with several numerical examples.

\subsection{High-Order Finite Element Method}\label{sec:hpfem}

In constrast with the standard finite element method ($h$-version of FEM) the high-order
finite element method adds a refinement parameter, the local polynomial order $p$,
hence the name $p$-version. When both refinements are available we refer to 
$hp$-version.
High-oder finite element methods have the capability for exponential convergence
provided the discretization is constructed properly in both domain (in $h$)
and local polynomial order (in $p$).

In this paper in all cases it is implicitly assumed that the exact parameterization of
the boundaries on the parameter space is known. This allows us to benefit
from efficient handling of large elements within the $p$-version
without significant loss of accuracy, and more importantly,
geometric refinements can be carried out with relative ease. 
This means that the number of 
elements can be kept relatively low.

Let us consider the Dirichlet problem \eqref{def_condensercap} and its 
weak solution $u_0$. 
The following theorem due to Babu{\v{s}}ka and Guo \cite{BaGuo},
sets the limit to the rate of convergence of the $hp$-FEM. 
Notice that construction of
the appropriate spaces is technical, but can be extended to parameterised surfaces. 
For rigorous treatment of the theory involved,
see Schwab \cite{schwab} and references therein. 
\begin{theorem} \label{propermesh}
Let the computational domain 
$G \subset \mathbb{R}^2$, $v$ the FEM-solution of \eqref{def_condensercap}, and
let the weak solution $u_0$ be in a suitable countably normed space where
the derivatives of arbitrarily high order are controlled.
Then
\[
\inf_v \|u_0 - v\|_{H^1(G)} \leq C\,\exp(-b \sqrt[3]{N}),
\]
where $C$ and $b$ are independent of $N$, the number of degrees of freedom. 
Here $v$ is computed on a proper geometric mesh, where the order of an individual
element is set to be its element graph distance to the nearest singularity.
(The result also holds for meshes with constant polynomial degree.)
\end{theorem}

There are many efficient error estimators available for $hp$-FEM.
The so-called auxiliary subspace error estimation fits particularly well
within our implementation.
Let $\cT$ be some $hp$-discretization on the
computational domain $G$. Assuming that the exact
solution $u \in H_0^1(G)$, defined on $\cT$, has finite energy, the
approximation problem is as follows: Find $\hat{u} \in V$ such that
\begin{equation}\label{eq:approximation}
  a(\hat{u},v)=l(v)\ (= a(u,v))\quad (\forall v \in V),
\end{equation}where
$a(\,\cdot\,,\,\cdot\,)$ and $l(\,\cdot\,)$, are the bilinear form
and the load potential, respectively. Additional degrees of
freedom are introduced by enriching the space $V$
via introduction of an auxiliary subspace or ``error
space'' $W \subset H_0^1(G)$ such that $V \cap W = \{0\}$. 
The error problem becomes thus: Find $\varepsilon \in W$ such that
\begin{equation}\label{eq:error}
  a(\varepsilon,v)=l(v)- a(\hat{u},v) (= a(u-\hat{u},v))\quad (\forall v \in W).
\end{equation}
This can be interpreted as a projection of the residual to the auxiliary space.

The main
result on this kind of estimators for the Dirichlet problem \eqref{def_condensercap}
is given in \cite{hno}. It should be mentioned that it is generally accepted that
the constant $K$ appearing in the theorem bwlow is in fact independent of $p$, yet
no rigorous proofs exist to support this observation.
\begin{theorem}[\cite{hno}]\label{KeyErrorThm}
There exist a constant $K$ depending only on the dimension $d$ and
polynomial degree $p$, continuity and coercivity constants $C$ and
$c$, and the shape-regularity of the triangulation $\cT$ such that
\begin{align*}
\frac{c}{C}\,\|\eE\|_1\leq\|u-\hat{u}\|_{1}\leq K(d,p)
\left(\|\eE\|_{1}+\osc(R,r,\cT)\right),
\end{align*}
where the residual oscillation depends on the volumetric and face
residuals $R$ and $r$, and the triangulation $\cT$.
\end{theorem}

\subsection{A Priori Refinement Strategies}\label{sec:refine}
In the classes of geometries considered here, two implementation challenges have to be met
to obtain the optimal exponential convergence.
First, for the slits, we must have geometric refinements of the elements at the singularities. Second, for
circular boundaries, 
the element boundaries must retain the exact parameterization.

\begin{figure}
\centering
 \subfloat{\includegraphics[height=2.25in]{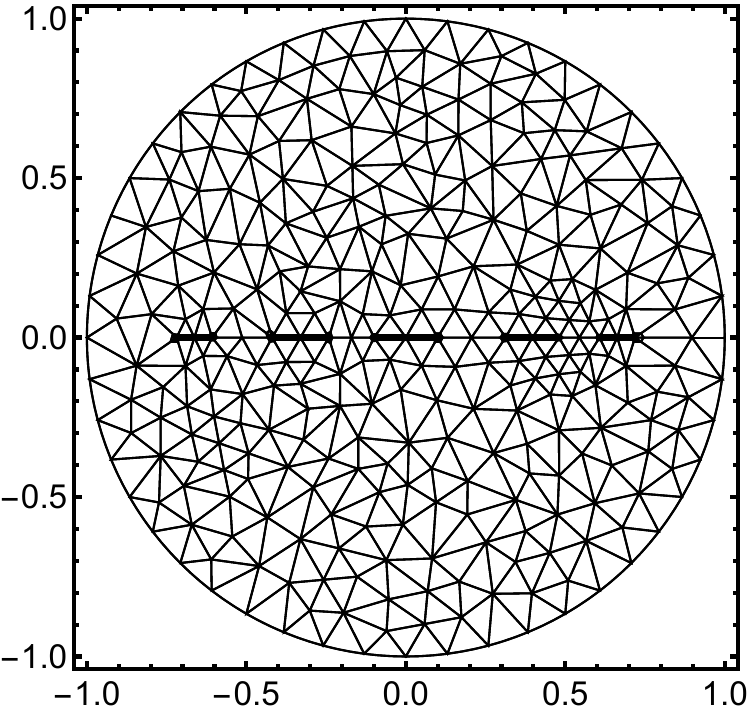}}
 \subfloat{\includegraphics[height=2.25in]{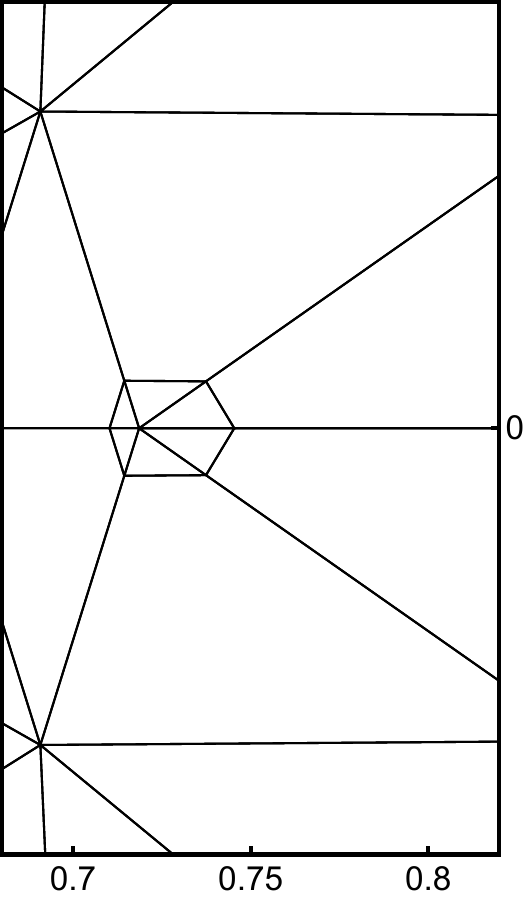}}
 \subfloat{\includegraphics[height=2.25in]{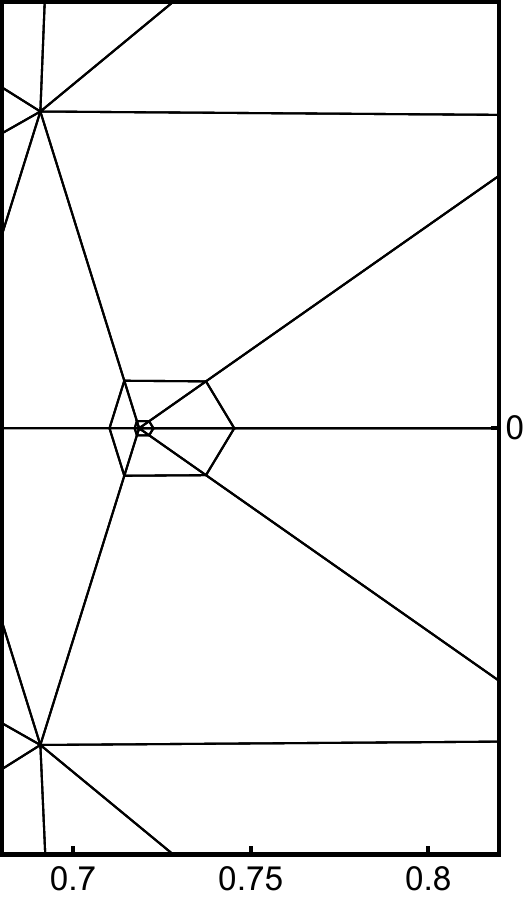}}
 \caption{Mesh refinements. Left: FEM mesh with the five segments
 on the diameter indicated with thick lines. Middle: A detail of the mesh
 at one of the end points of the segments after one application of the
 replacement rule. Right: Successive levels of refinements, eight altogether,
 are shown in the plot (the smallest ones are not visible in the given scale).
  }\label{fig:refinement}
\end{figure}

In problems with singularities with known locations, 
a priori optimally refined meshes can be computed
using rule based algorithms \cite{HaTu}. 
The geometric refinement process is illustrated in
Figure~\ref{fig:refinement}: First a geometrically conforming mesh is
generated such that every singularity is isolated, that is, 
local refinement rules can be applied simultaneously without
violating the conformity of the mesh.
Each local refinement can be applied multiple times in a nested fashion.
One of the drawbacks of this approach is that 
unwinding local refinements is difficult, since
typical geometric invariants of the triangulations are not valid within the local
refinements, for instance, the Delaunay property (maximization of the minimal angle).
In the solution process we are content to adapt the discretization simply by modifying
the a priori strategy, in other words, by remeshing the whole domain.

\begin{figure}
\centering
 \subfloat{\includegraphics[height=2.25in]{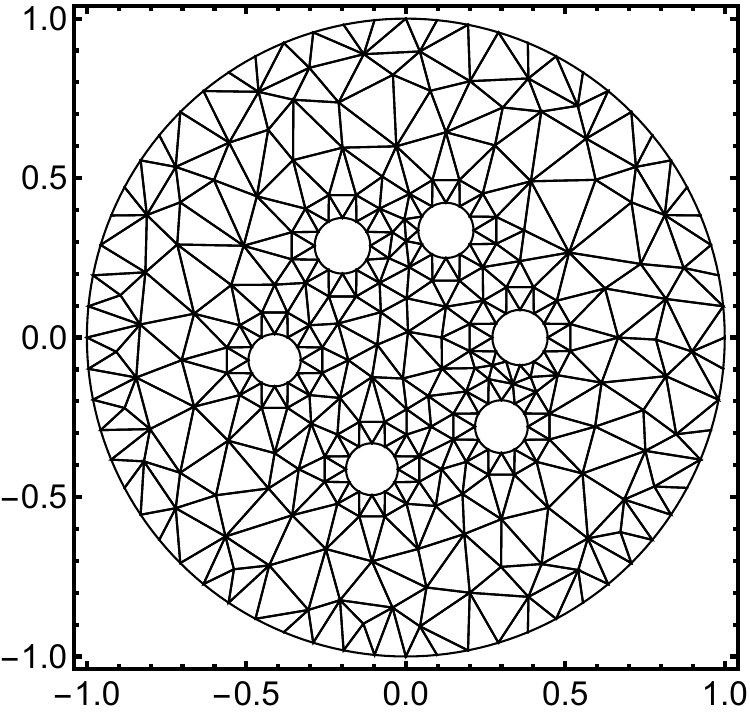}}\qquad
 \subfloat{\includegraphics[height=2.25in]{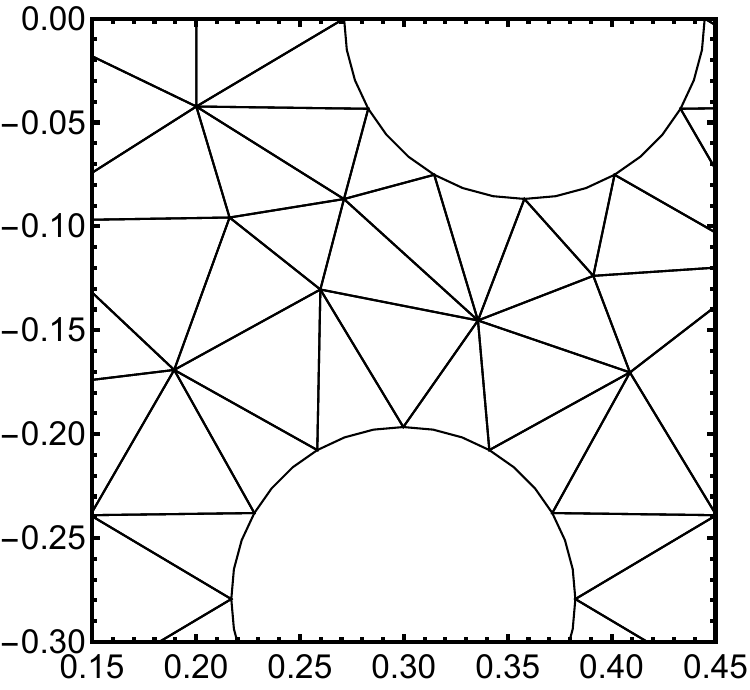}}
 \caption{ Curved elements. Left: FEM mesh with six disks. Right:
A detail showing the curved edges.  }\label{fig:refinementcurved}
\end{figure}

For problems with circular boundaries, each element at the boundary
is curved with geometrically faithful parameterization.
One example with associated detail is shown in 
Figure~\ref{fig:refinementcurved}.

\subsection{BIE method}\label{sec:bie}
We shall consider two types of condensers $(\B,E)$ in this paper. In the first type, the compact set $E$ is assumed to be the union of $m$ disjoint hyperbolic disks in $\B$. For the second type, we assumed that $E$ is the union of $m$ disjoint hyperbolic segments in $\B$. 
For both types of domains, the capacity $\capa(\B,E)$ can be computed using the boundary integral equation (BIE) method presented in~\cite{nv-g}. The method is based on the BIE with the generalized Neumann kernel. This method is briefly reviewed in this section.
However, before implementing the numerical method, we first convert the hyperbolic disks and segments to Euclidean ones, using~\eqref{hkv4.20}.

\subsubsection{Domains bounded by smooth curves}\label{sec:smooth}

When $E$ is a union of $m$ disjoint hyperbolic disks, then the domain $\Omega=\B\backslash E$ is a bounded multiply connected domain of connectivity $m+1$ whose boundaries are circles. The orientation of the external circle $C_0$ is counterclockwise oriented and the inner circles $C_1,\ldots,C_m$ are clockwise oriented. The external circle $C_0$ is parametrized by $\eta_0(t)$ for $t\in J_0=[0,2\pi]$. Each inner circle $C_j$ is parametrized by $\eta_j(t)$, $t\in J_j=[0,2\pi]$, for $j=1,2,\ldots,m$. 
Let $J$ be the disjoint union of the $m+1$ intervals $J_j=[0,2\pi]$, $j=0,1,\ldots,m$. We define a parameterization of the whole boundary $C=\cup_{j=0}^m C_j$ on $J$ by (see~\cite{Nas-ETNA} for the details)
\[
\eta(t)=\left\{
\begin{array}{cc} 
	\eta_0(t), & t\in J_0, \\ 
	\eta_1(t), & t\in J_1, \\
	\vdots \\
	\eta_m(t), & t\in J_m. \\ 
\end{array}
\right.
\]

With the parameterization $\eta(t)$ of the whole boundary $C$, we define a complex function $A$ by
\begin{equation}\label{eq:A}
	A(t) = \eta(t)-\alpha,
\end{equation}
where $\alpha$ is a given point in the domain $\Omega$. 
For each $k=1,2,\ldots,m$, let $z_k$ be a given point interior to the circle $C_k$, let the function $\gamma_k$ be defined by
\begin{equation}\label{eq:gam-k}
	\gamma_k(t)=\log|\eta(t)-z_k|,
\end{equation}
and let $\mu_k$ be the unique solution of the BIE
\begin{equation}\label{eq:ie}
	\mu_k-\bN\mu_k=-\bM\gamma_k,
\end{equation}
where $\bN$ is the integral operator with the generalized Neumann kernel 
\begin{equation}\label{eq:N}
	N(s,t) :=
	\frac{1}{\pi}\Im\left(\frac{A(s)}{A(t)}\frac{\eta'(t)}{\eta(t)-\eta(s)}\right), \quad (s,t)\in J\times J,
\end{equation}
and $\bM$ is the integral operator with the kernel
\begin{equation}\label{eq:M}
	M(s,t) :=
	\frac{1}{\pi}\Re\left(\frac{A(s)}{A(t)}\frac{\eta'(t)}{\eta(t)-\eta(s)}\right),
	\quad (s,t)\in J\times J.
\end{equation}
Then the function $h_k$ given by
\begin{equation}\label{eq:h}
	h_k=[\bM\mu_k-(\bI-\bN)\gamma_k]/2
\end{equation}
is a piecewise constant function, i.e., 
\[
h_k(t)=\left\{
\begin{array}{cc} 
	h_{0,k}, & t\in J_0, \\ 
	h_{1,k}, & t\in J_1, \\
	\vdots \\
	h_{m,k}, & t\in J_m, \\ 
\end{array}
\right.
\]
where $h_{0,k},h_{1,k},\ldots,h_{m,k}$, $k=1,2,\ldots,m$, are real constants. The capacity $\capa(\B,E)$ can be then computed by~\cite[Eq.~(3.9)]{nv-g} 
\begin{equation}\label{eq:cap-ak}
	\capa(\B,E)=2\pi\sum_{k=1}^m a_k,
\end{equation}
where the values of the $m$ real constants $a_1,\ldots,a_{m}$ are computed by solving the $(m+1)\times(m+1)$ linear system
\begin{equation}\label{eq:sys-method}
	\left[\begin{array}{ccccc}
		h_{0,1}    &h_{0,2}    &\cdots &h_{0,m}      &1       \\
		h_{1,1}    &h_{1,2}    &\cdots &h_{1,m}      &1       \\
		\vdots     &\vdots     &\ddots &\vdots       &\vdots  \\
		h_{m,1}    &h_{m,2}    &\cdots &h_{m,m}      &1       \\
	\end{array}\right]
	\left[\begin{array}{c}
		a_1    \\a_2    \\ \vdots \\ a_{m} \\  c 
	\end{array}\right]
	= \left[\begin{array}{c}
		0 \\  1 \\  \vdots \\ 1  
	\end{array}\right].
\end{equation}
The constants $b_1,\ldots,b_m$ in~\eqref{eq:cap-bk} are related to the constants $a_1,\ldots,a_m$ by
\begin{equation}\label{eq:bk-ak}
b_k=2\pi a_k, \quad k=1,2,\ldots,m.
\end{equation}

The BIE~\eqref{eq:ie} can be discretized by the Nystr\"om method with the trapezoidal rule to obtain an $(m+1)n\times(m+1)n$ linear system where $n$ is the number of the discretization points in each boundary component. The linear system can then be solved by the MATLAB function {\tt gmres} and the matrix-vector product in {\tt gmres} can be computed in $O(mn)$ operations using the MATLAB function {\tt zfmm2dpart} from the fast multipole method (FMM) MATLAB toolbox FMMLIB2D~\cite{Gre-Gim12}. 
The boundary components of the domain $\Omega$ are circles and hence the integrands in~\eqref{eq:ie} and~\eqref{eq:h} are analytic. Thus the trapezoidal rule converges exponentially with $O\left(e^{-\sigma n}\right)$ when it is used to discretize the integrals in~\eqref{eq:ie} and~\eqref{eq:h}~\cite{Tre}. The constant $\sigma$ depends on the configuration of the domain $\Omega$ as well as the auxiliary point $\alpha$.
The numerical solution of the integral equation will then converge with a similar rate of convergence as the trapezoidal rule~\cite{Atk}. 
This method for solving the BIE~\eqref{eq:ie} was implemented in the MATLAB function {\tt fbie} presented in~\cite{Nas-ETNA}. 
The MATLAB function {\tt fbie} provides us with approximations to the solution $\mu_k$ of the BIE~\eqref{eq:ie} as well as the piecewise constant function $h_k$ in~\eqref{eq:h}. The computed values of $h_k$ are used to set up the $(m+1)\times(m+1)$ linear system~\eqref{eq:sys-method}, which will be solved using the Gauss elimination method (here $m+1$ is the number of boundary components of the domain $G$ which is usually small). 
By computing the constants $a_1,\ldots,a_m$, the value of the capacity $\capa(G,E)$ is given by~\eqref{eq:cap-ak}. Further, the values of the constants $b_1,\ldots,b_m$ are given by~\eqref{eq:bk-ak}.
See~\cite{Nas-ETNA,nv-g} for details.

\subsubsection{Domains bounded by slits}\label{sec:slits}

The BIE method presented above can be used to compute the capacity of only condensers bounded by smooth or piecewise Jordan curves~\cite{Nas-ETNA,nv-g}. Since the Dirichlet integral is conformally invariant, the capacities for the cases for which the plates of the condenser are slits can be computed with the help of conformal mappings. In this paper, we consider two types of domains $\Omega$ bounded by slits. 

\begin{figure}
	\centering
{
		\includegraphics[width=0.4\textwidth]{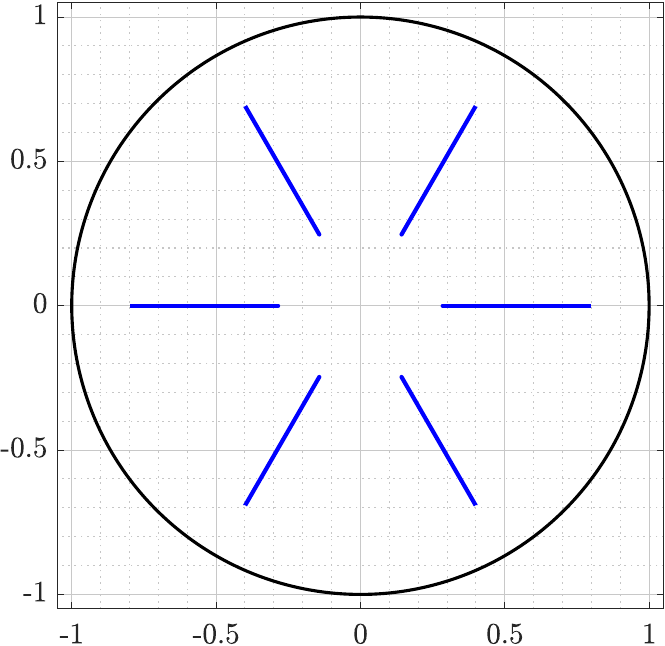}
\hfill
		\includegraphics[width=0.4\textwidth]{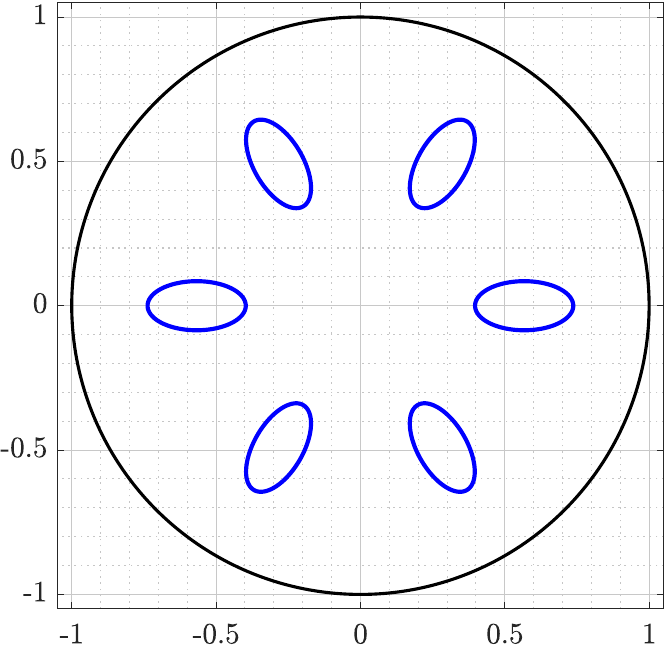}
	}
	\caption{A schematic of a given bounded multiply connected domain $\Omega$ interior to the unit circle and exterior to $m$ radial slits (left) and its conformally equivalent computed domain $D$ interior to the unit circle and exterior to $m$ smooth Jordan curves (right) for $m=6$.}\label{fig:conf-rad}
\end{figure}

In the first case, we assume that $\Omega$ is the unit disk with $m$ radial slits. For such a case, we can use the iterative method presented in~\cite{NG18} to compute a conformally equivalent domain $D$ bounded by smooth Jordan curves so that our method presented in Section~\ref{sec:smooth} can be used. A schematic of the domain $\Omega$ and its conformally equivalent computed domain $D$ for $m=6$ is presented in Figure~\ref{fig:conf-rad}. 

\begin{figure}
	\centering
{
		\includegraphics[width=0.4\textwidth]{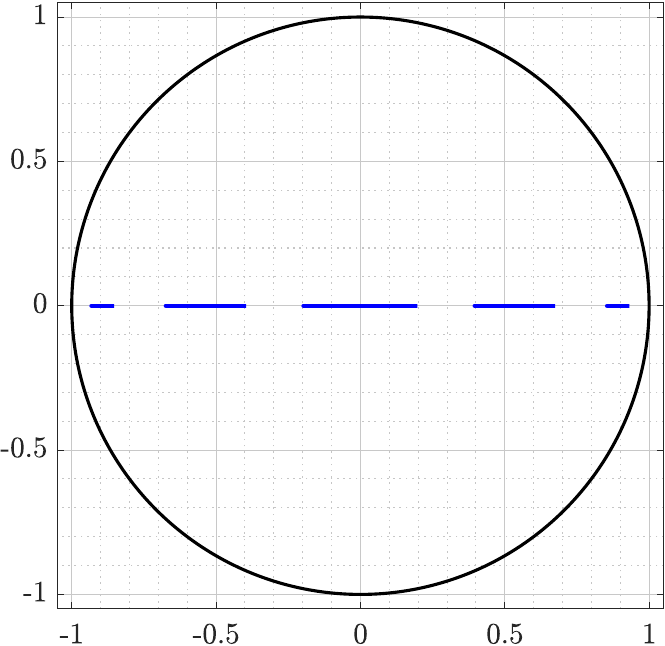}
\hfill
		\includegraphics[width=0.4\textwidth]{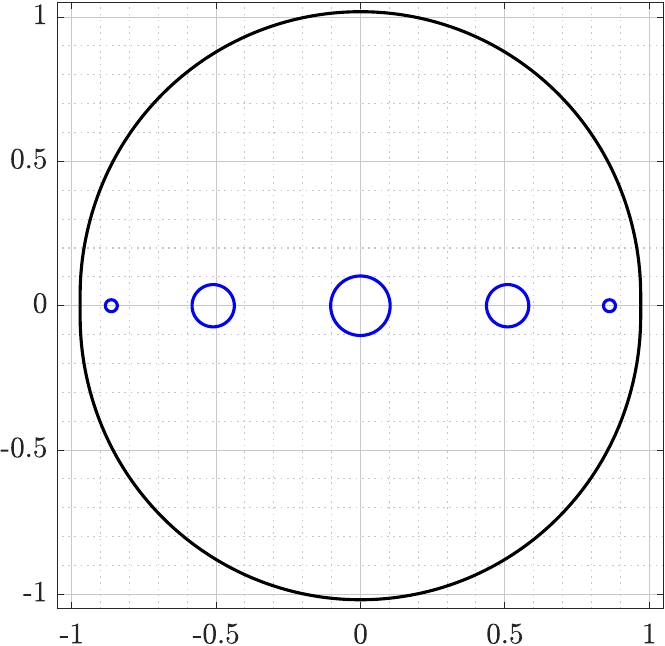}
	}
	\caption{A schematic of a given bounded multiply connected domain $\Omega$ interior to the unit circle and exterior to $m$ rectilinear slits (left) and its conformally equivalent computed domain $D$ bounded by $m+1$ smooth Jordan curves (right) for $m=5$.}\label{fig:conf-rec}
\end{figure}

\begin{figure}
	\centering
	{
		\includegraphics[width=0.4\textwidth]{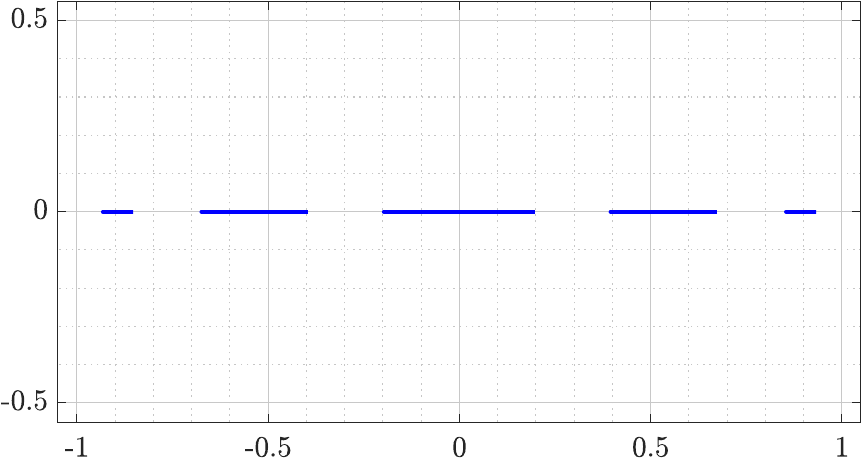}
\hfill
		\includegraphics[width=0.4\textwidth]{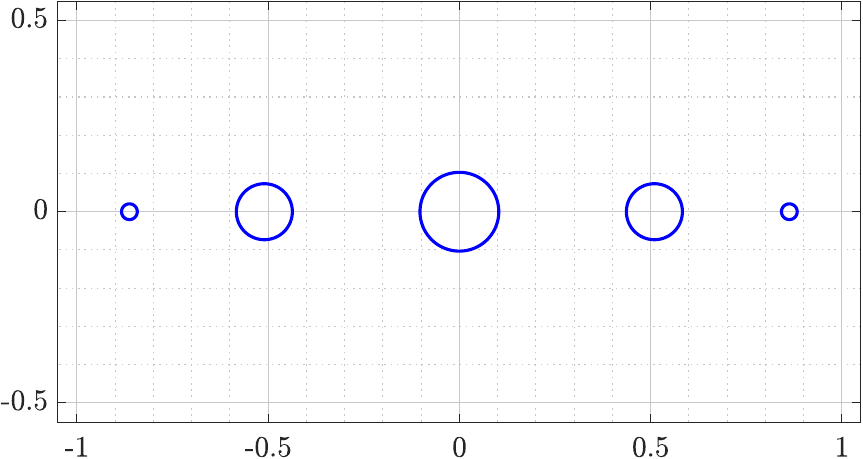}
	}
	\caption{A schematic of a given unbounded multiply connected domain $\hat\Omega$ exterior to $m$ rectilinear slits (left) and its conformally equivalent computed domain $\hat D$ exterior to $m$ smooth Jordan curves (right) for $m=5$.}\label{fig:conf-recu}
\end{figure}

In the second case,  we assume that $\Omega$ is the unit disk with $m$ rectilinear slits on the real line (see Figure~\ref{fig:conf-rec} (left) for $m=5$). Unlike the domain in the first case, this domain is not one of the canonical slit domains (see~\cite{Koe,Nas-JMAA11}). Thus, in this case, we first consider the unbounded domain $\hat\Omega$ in the exterior of the $m$ rectilinear slits which is a canonical slit domain (see Figure~\ref{fig:conf-recu} (left) for $m=5$). We use the iterative method presented in~\cite{NG18} to compute a conformally equivalent domain $\hat D$ in the exterior of $m$ smooth Jordan curves and the conformal mapping $w=\Phi(z)$ from the domain $\hat D$ onto $\hat\Omega$  (see Figure~\ref{fig:conf-recu} (right) for $m=5$). Hence, $z=\Phi^{-1}(w)$ is a conformal mapping from the domain $\hat\Omega$ onto $\hat D$. Since the unit circle is in the interior of the domain $\hat\Omega$, the conformal mapping $z=\Phi^{-1}(w)$ can be used to compute the image of the unit circle which will be a smooth Jordan curve exterior to the computed $m$ smooth Jordan curves. Thus, the conformal mapping $z=\Phi^{-1}(w)$ maps the given domain $\Omega$ onto a conformally equivalent domain $D$ bounded by smooth Jordan curves so that the method reviewed in Section~\ref{sec:smooth} can be used  (see Figure~\ref{fig:conf-rec} (right) for $m=5$). Notice that the external curve in Figure~\ref{fig:conf-rec} (right) is not a circle.

For details on the iterative method for computing the domain $D$ for both cases of slit domains discussed above, we refer the reader to~\cite{NG18} (see also~\cite{knv} for other types of slit domains).

\subsection{Nonlinear Optimization: Interior-Point Method}\label{sec:ipm}

The numerical optimization algorithm of our choice is the interior-point method
as implemented in Mathematica (\texttt{FindMaximum}, \cite{Wol}) and Matlab (\texttt{fmaxcon},  \cite{Matlab}).
The task is to find an optimal configuration for a constellation of 
hyperbolic disks $E$ with fixed radii, where at every step the current configuration
is solved using either one of the methods described above.
The standard textbook reference is Nocedal and Wright \cite{Noc}.

\begin{figure}
	\centering
	{
		\includegraphics[width=0.5\textwidth]{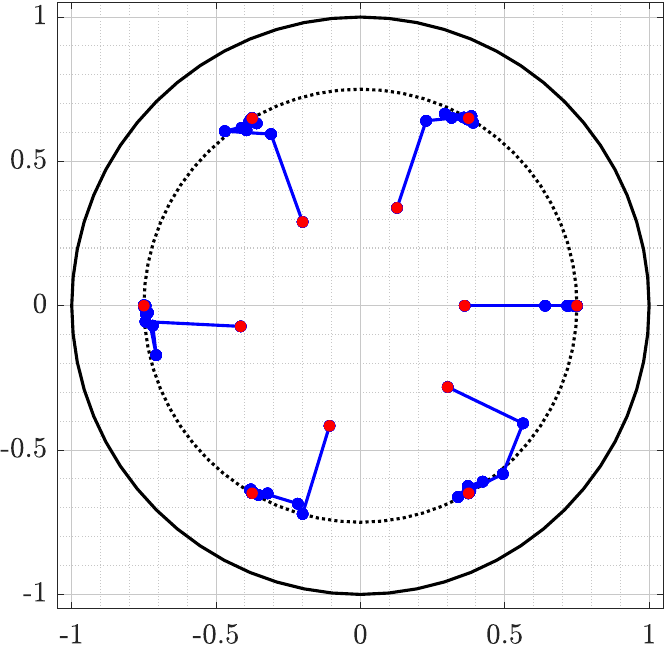}
	}\hfill
	\caption{The positions of the centers of the six disks in Figure~\ref{fig:6Dint} during the optimization process using the BIE method. The red dots corresponding to the initial positions (Figure~\ref{fig:6Dint} (left)) and the maximal positions  (Figure~\ref{fig:6Dint} (right)).}\label{fig:path}
\end{figure}

In the most general case the problem is defined as in \eqref{eq:max},
where the only constraints are geometric ones, that is, the
disks are not allowed to overlap, and they are not allowed to drift to the boundary, for instance,
they must lie within a disk $B^2(0,R)$ with same prescribed radius $R$, or alternatively their centers must 
lie inside some constraining disk. The radii are fixed and the
optimization concerns only the locations of the disks.
For example, for the constellation of six hyperbolic disks with the initial positions shown in Figure~\ref{fig:6Dint} (left), the positions of the centers of the six disks during the optimization process are shown in Figure~\ref{fig:path} where the red dots correspond to the initial positions (Figure~\ref{fig:6Dint} (left)) and the maximal positions (Figure~\ref{fig:6Dint} (right)).

The maximization problem is formally defined as
\begin{alignat}{3}
\max_E              &\quad&  {\rm cap}(\mathbb{B}^2,E)\quad  &&    & \nonumber \\
\text{subject to: } &\quad&  E_i \cap E_j  = \emptyset\quad  &&    &\quad \forall\ i,j=1,\ldots,m, i \neq j \label{eq:max}\\
                    &\quad&  E_j  \subset B^2(0,R)      &&    &\quad \forall\ j=1,\ldots,m. \nonumber 
\end{alignat}
This nonlinear optimization problem can be solved using the interior-point method, and 
the solution would be a local maximum.

Notice, that the objective function is indeed the capacity of the
constellation. The number of evaluations needed for the optimization is greater than the number of iteration
steps, since the gradients and Hessians must be approximated numerically.
One of the insights gained over many such computations is that
the optimization depends on the high accuracy of
the capacity solver, since otherwise the approximate derivatives are not sufficiently accurate.

In the context of this work, there have been no attempts to devise a special
method that would incorporate some of the insights gathered during this study.
Instead, the numerical optimization is used to challenge those insights and
therefore the optimizations have been computed with minimal input information.

\section{Numerical experiments: Constellations of circular domains}

In this section the focus is on constellations of disks. In the maximization of the capacity the positions of the disks
are subject to two types of geometric constraints, they are either constrained to a disk of given radius or an interval of fixed length
on the real line. The experiments in turn either cover full parameter ranges or are general in the sense that the initial configurations 
are random, but satisfy the constraints, of course.
We first consider constellations of two disks of equal hyperbolic radii, and then extend the investigation to constellations with six disks constrained to a disk, and to constellations with five disks with centers constrained to an interval. In the two latter cases also the case of unequal hyperbolic radii is studied.
In the final experiment the constellation is condensed into a single disk with equal
capacity. The objective is to compare the hyperbolic area and perimeter of a constellation to that of a condensed one.

\subsection{Constellation of two disks with constrained positions}

\begin{figure}
	\centering
	{
		\includegraphics[width=0.4\textwidth]{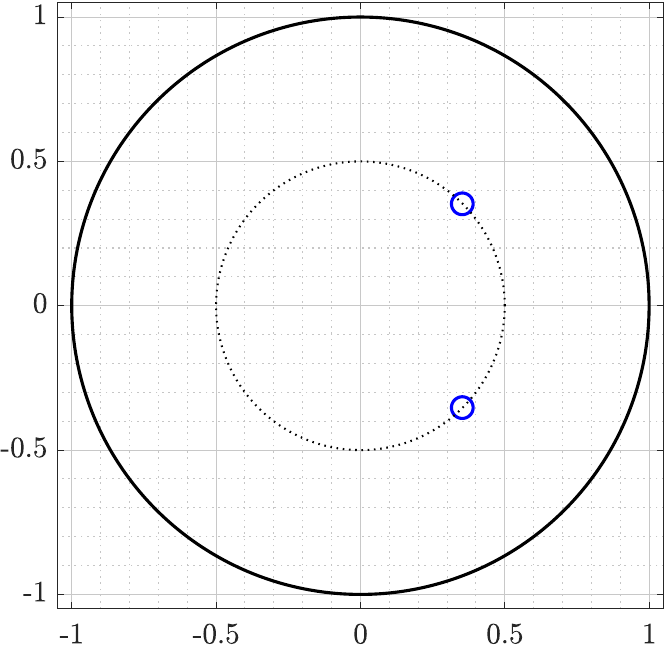}
	}\hfill
	{
		\includegraphics[width=0.4\textwidth]{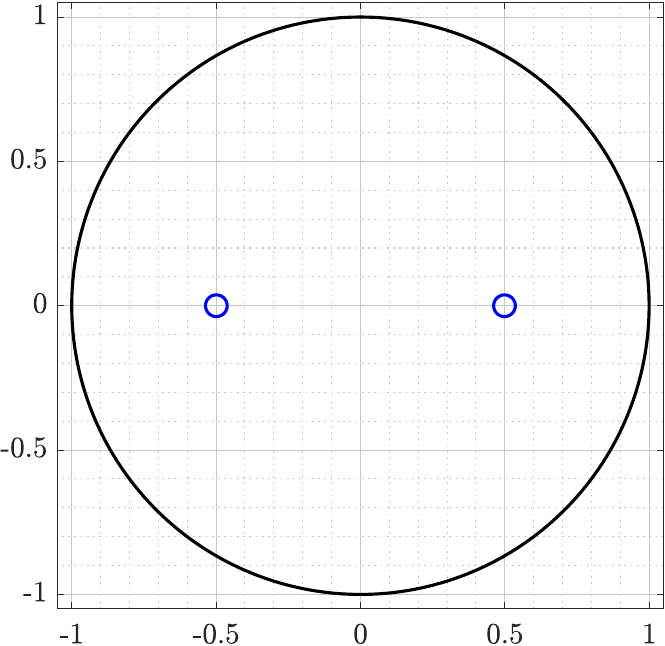}
	}
	\caption{Two hyperbolic disks with hyperbolic radius $0.1$. The centers of the two disks are $Re^{\pm\i\theta}$ for $R=0.5$ and $\theta=\pi/4$ (left) and $\pm x$ for $x=0.5$ (right).}\label{fig:D2Capfig}
\end{figure}

\begin{figure}
	\centering
	{
		\includegraphics[width=0.4\textwidth]{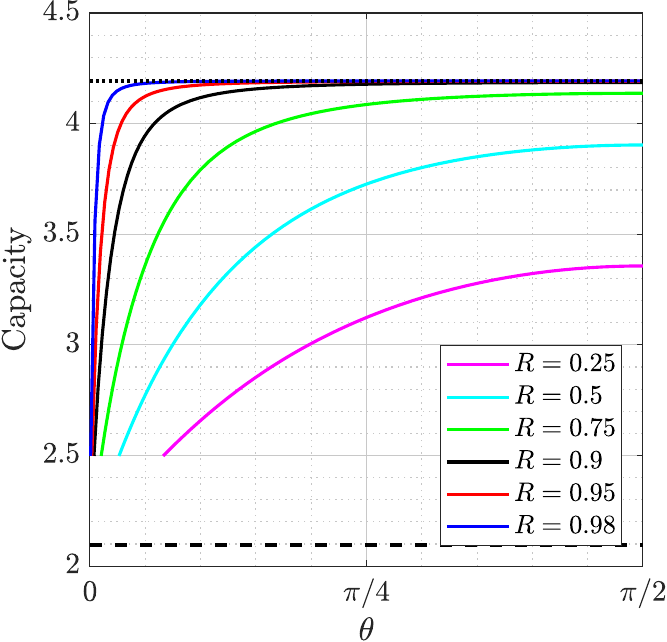}
	}\hfill
	{
		\includegraphics[width=0.4\textwidth]{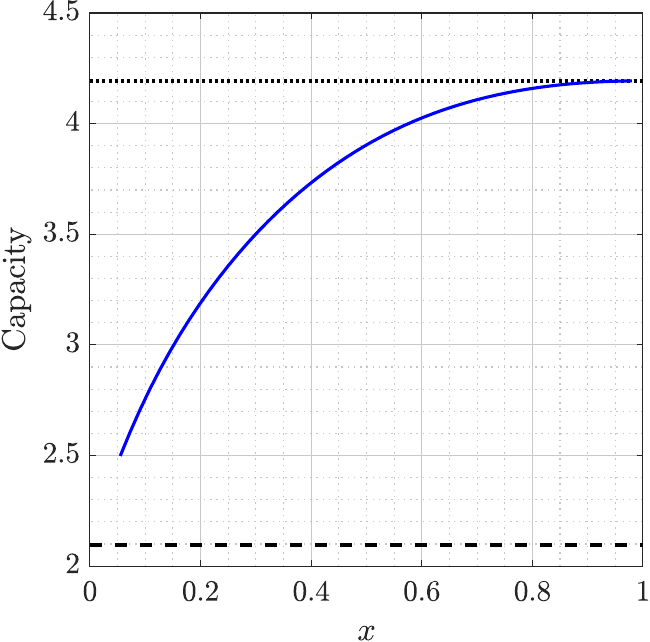}
	}
	\caption{The capacity for two hyperbolic disks with hyperbolic radius $r = 0.1$ . 
	The centers of the two disks are $Re^{\pm\i\theta}$ (left) and $\pm x$ (right),
	where the leftmost point on the curve corresponds to Cartesian coordinate
	$\th(r/2) \approx 0.05$, that is, when the two disks touch.
	}\label{fig:D2Cap}
\end{figure}

We begin with the constellation $E$, union of two hyperbolic disks $D_1$ and $D_2$ with equal hyperbolic radius $r$. 
First we assume that the centers of these disks are on $Re^{\pm\i\theta}$ where 
$0<\theta_{\rm min}<\theta<\pi-\theta_{\rm min}<\pi$ and
\[
\theta_{\rm min} = \arcsin\left(\frac{(1-R^2)\sh r}{2R}\right).
\] 
See Figure~\ref{fig:D2Capfig} (left) for $R=0.5$ and $\theta=\pi/4$.
The two disks touch each other when $\theta=\theta_{\rm min}$ or $\theta=\pi-\theta_{\rm min}$. When $r=0.1$, the values of $\capa(\mathbb{B}^2,E)$ vs. $\theta$ are shown in Figure~\ref{fig:D2Cap} (left) for several values of $R$. 

Then we assume that the centers of these disks are on $\pm x$ where $x_{\rm min}<x<1$ and $x_{\rm min}=\th(r/2)$ where the two disks touch each other when $x=x_{\rm min}$ (See Figure~\ref{fig:D2Capfig} (right) for $x=0.5$). The values of $\capa(\mathbb{B}^2,E)$ for $r=0.1$ vs. $x$ are shown in Figure~\ref{fig:D2Cap} (right).

Note that $\capa(\mathbb{B}^2,D_i)=-2\pi/\log(\th(r/2))$, $i=1,2$. For $r=0.1$, the values of $\capa(\mathbb{B}^2,D_1)$ are shown in Figure~\ref{fig:D2Cap} as ``dashed line'' and the values of  $\capa(\mathbb{B}^2,D_1)+\capa(\mathbb{B}^2,D_2)$ as ``dotted line.''

\subsection{Constellation of six disks constrained to a disk}\label{sec:6D}

We increase the number of disks and 
consider the positions of a constellation of six hyperbolic disks that maximize the
capacity $\capa(\mathbb{B}^2,E)$ under the constraint that the hyperbolic centers
of these disks are inside the Euclidean disk $|z|\le R$ (we assume in the
examples below that $R=0.75$). The disks are numbered $D_1$ to $D_6$ in
counterclockwise orientation. We denote the center of the disk $D_j$ by
$z_j$, $j=1,\ldots,6$.  Without any loss of generality, we assume that the
center $z_1$ of the disk $D_1$ lies on the positive real axis. 

\begin{figure}
	\centering
	{
		\includegraphics[width=0.31\textwidth]{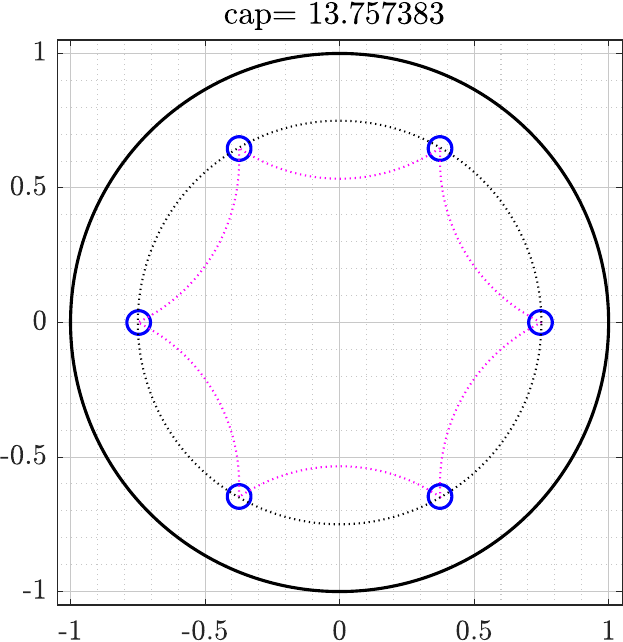}
	}
	{
		\includegraphics[width=0.31\textwidth]{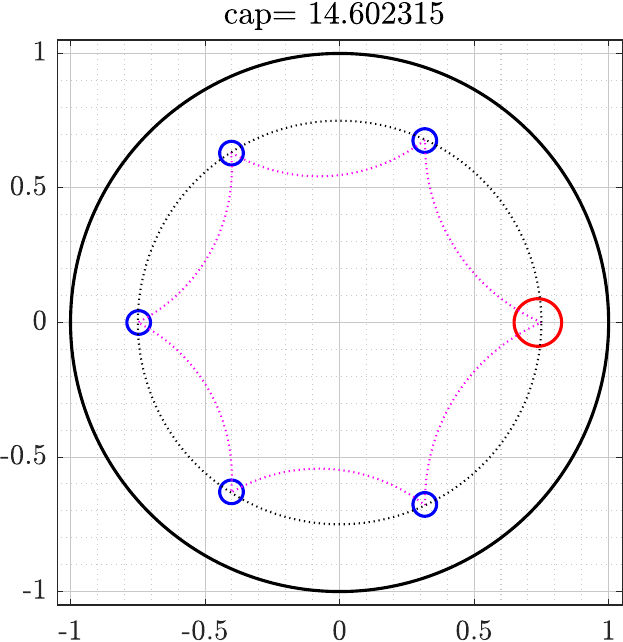}
	}
	{
		\includegraphics[width=0.31\textwidth]{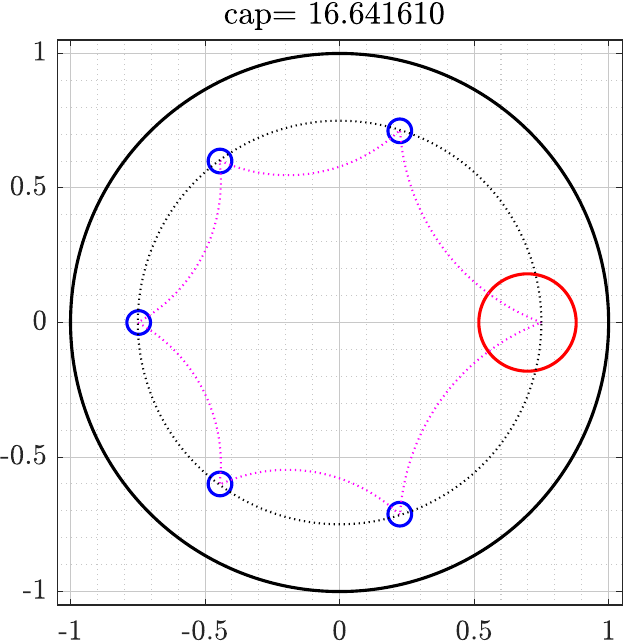}
	}
	\caption{Disk 1 has hyperbolic radius $0.2$ (left), $0.4$ (center), and $0.8$ (right). Disks 2 through 6 have hyperbolic radii 0.2.}\label{fig:6Dn1}
\end{figure}

First we assume that all six disks have equal hyperbolic radii $= 0.2$,
and the initial positions are random within the given constraints.
The configuration which maximizes the capacity $\capa(\mathbb{B}^2,E)$ 
has the maximal dispersion property:
The positions of these six disks are on the Euclidean circle $|z|=R$ 
and, moreover, are symmetric, that is, the
hyperbolic distances between the centers of any two adjacent disks are
equal (see Figure~\ref{fig:6Dn1} (left) and Table~\ref{tab:6Dn1}).
The computed capacity $= 13.757381$.

\begin{figure}
	\centering
	\includegraphics[width=0.4\textwidth]{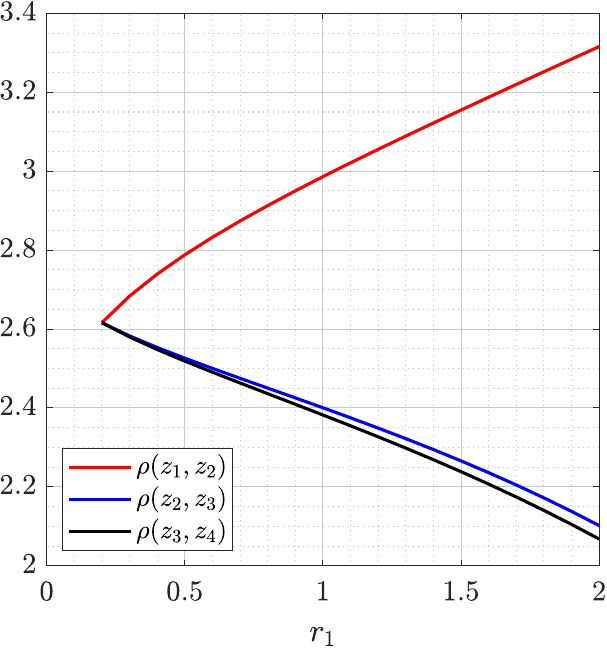}
	\hfill
	\includegraphics[width=0.4\textwidth]{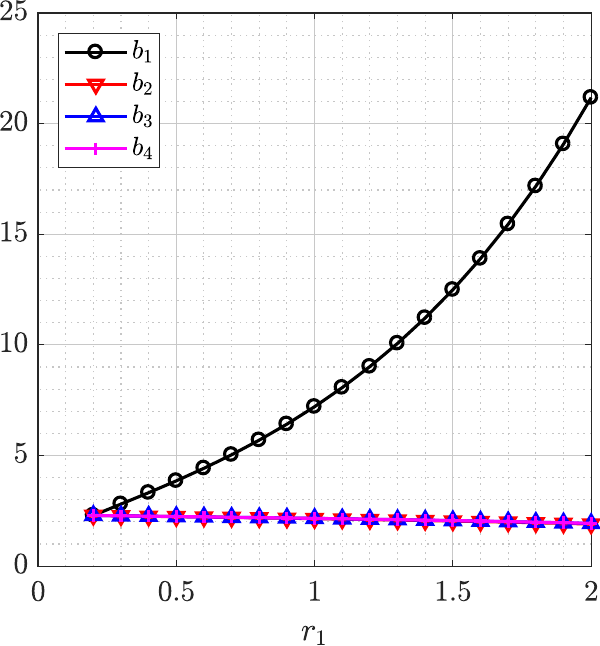}
	\caption{The hyperbolic distances between the centers of consecutive disks (left) and the values of the constants $b_1$, $b_2$, $b_3$, and $b_4$ in~\eqref{eq:cap-bk} (right) as functions of $r_1$.}\label{fig:6Dnv}
\end{figure}

When the hyperbolic radius of one of these disks is changed either to $0.4$ (see
Figure~\ref{fig:6Dn1} (center)) or $0.8$ (see Figure~\ref{fig:6Dn1} (right)),
the centers of the other disks move away from the larger disk (see
Table~\ref{tab:6Dn1}), yet rotational symmetry is preserved for the maximal configuration. 
To study closely the impact of increasing the hyperbolic radius of only one disk on the positions that maximize the capacity $\capa(\mathbb{B}^2,E)$, we assume that the hyperbolic radius of the first disk $D_1$ is $r_1$ and the hyperbolic radii of the remaining five disks $D_2$--$D_6$ are $0.2$. As above, we find the positions of these six disks that maximize the capacity $\capa(\mathbb{B}^2,E)$ under the above constraint. For the positions that maximizes the capacity $\capa(\mathbb{B}^2,E)$, we compute the hyperbolic distances $\rho(z_1,z_2)$, $\rho(z_2,z_3)$, and $\rho(z_3,z_4)$ and the values of the constants $b_1$, $b_2$, $b_3$, and $b_4$ in~\eqref{eq:bk-ak} where the values of $r_1$ are changing from $0.2$ to $2$. The obtained numerical results are presented in Figure~\ref{fig:6Dnv}. 
Notice that the constant $b_k$ can be regarded as the contribution of the disk set $D_k$ to the capacity ${\rm cap}(G,E)$, for $k=1,2,\ldots,m$. As we can see, the values of $b_1$ increased as $r_1$ increased and the values of $b_2$, $b_3$, and $b_4$ are almost constants.	
Notice also that, due to symmetry, $\rho(z_1,z_2)=\rho(z_1,z_6)$, $\rho(z_2,z_3)=\rho(z_6,z_5)$, $\rho(z_3,z_4)=\rho(z_5,z_4)$, $b_5=b_3$, and $b_6=b_2$.

\begin{figure}
	\centering
	{
		\includegraphics[width=0.31\textwidth]{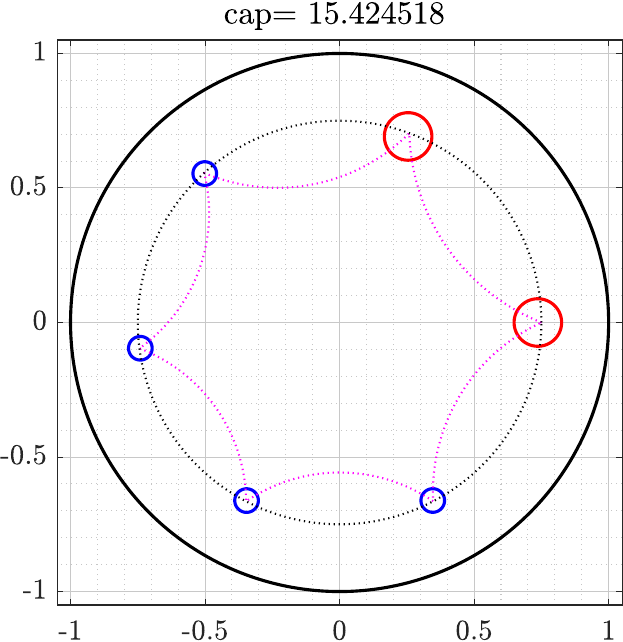}
	}
	{
		\includegraphics[width=0.31\textwidth]{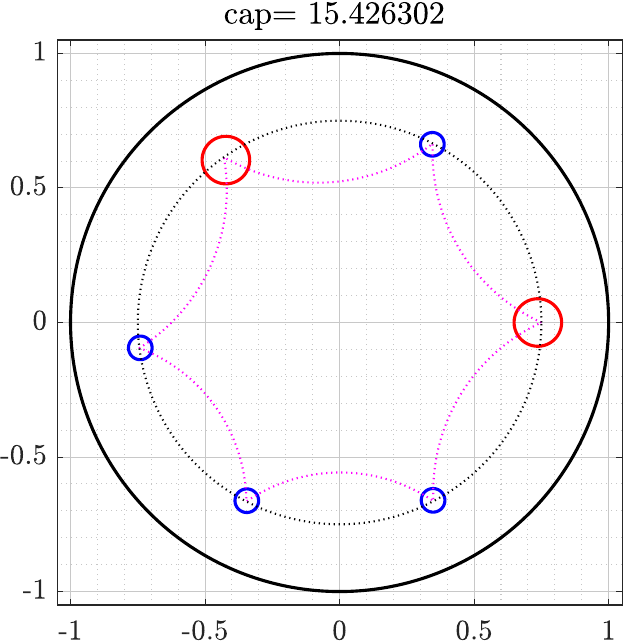}
	}
	{
		\includegraphics[width=0.31\textwidth]{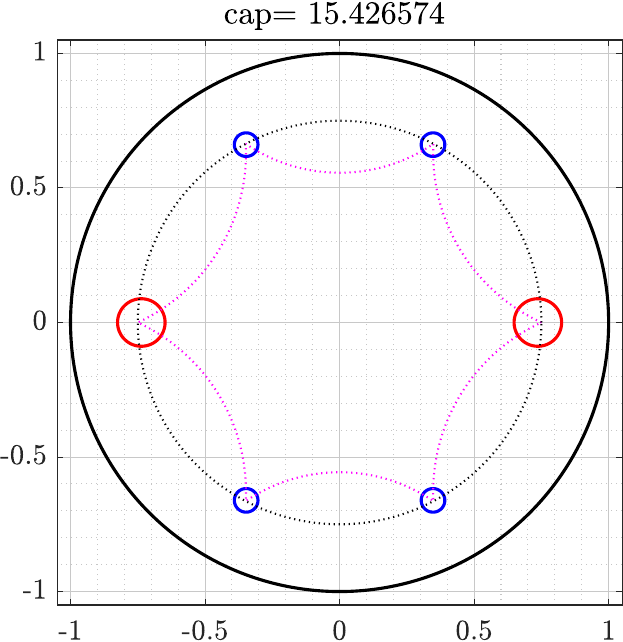}
	}\\
	{
		\includegraphics[width=0.31\textwidth]{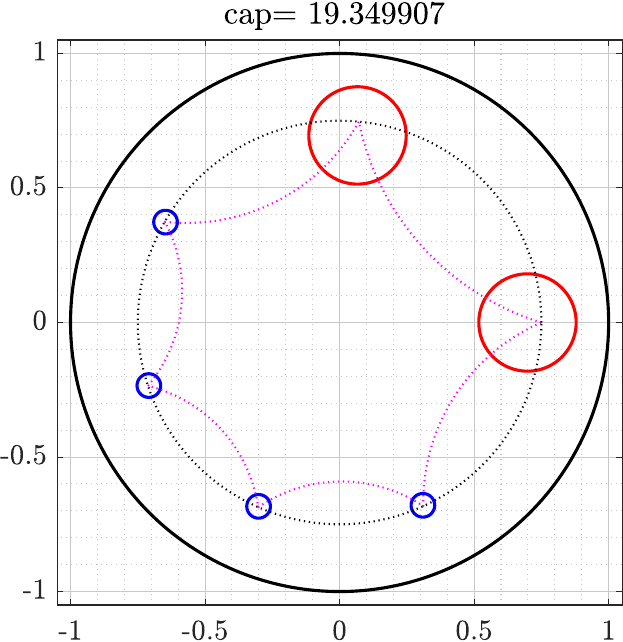}
	}
	{
		\includegraphics[width=0.31\textwidth]{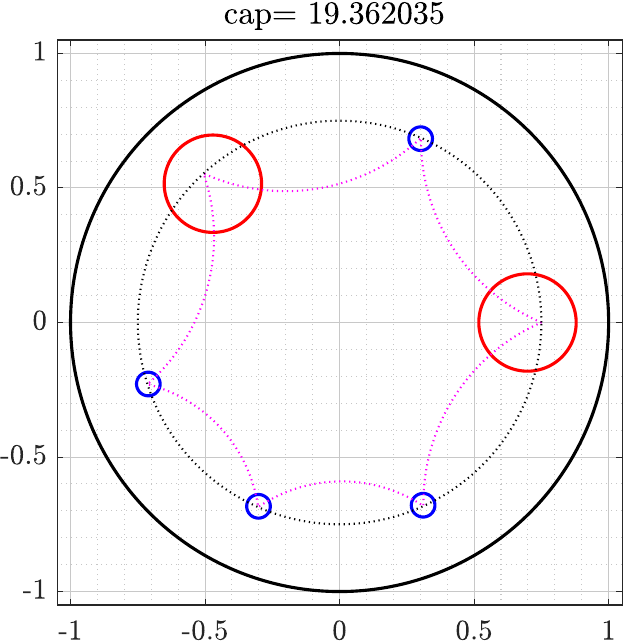}
	}
	{
		\includegraphics[width=0.31\textwidth]{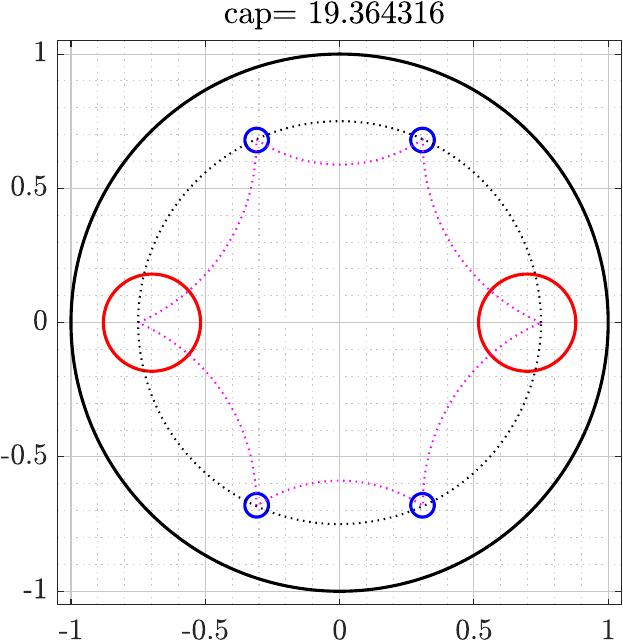} 
	}
	\caption{First row: Four disks with hyperbolic radii 0.2 and two disks with hyperbolic radii 0.4. Second row: Four disks with hyperbolic radii 0.2 and two disks with hyperbolic radii 0.8.}\label{fig:6Dn2}
\end{figure}

\begin{figure}
	\centering
	{
		\includegraphics[width=0.31\textwidth]{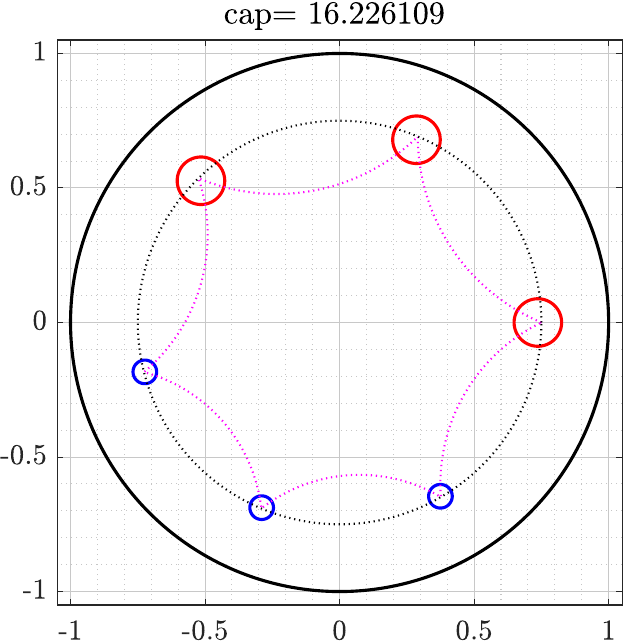}
	}
	{
		\includegraphics[width=0.31\textwidth]{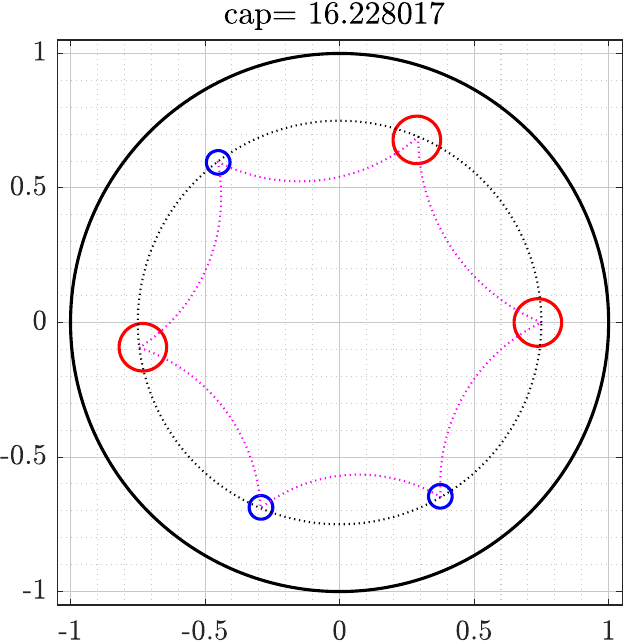}
	}
	{
		\includegraphics[width=0.31\textwidth]{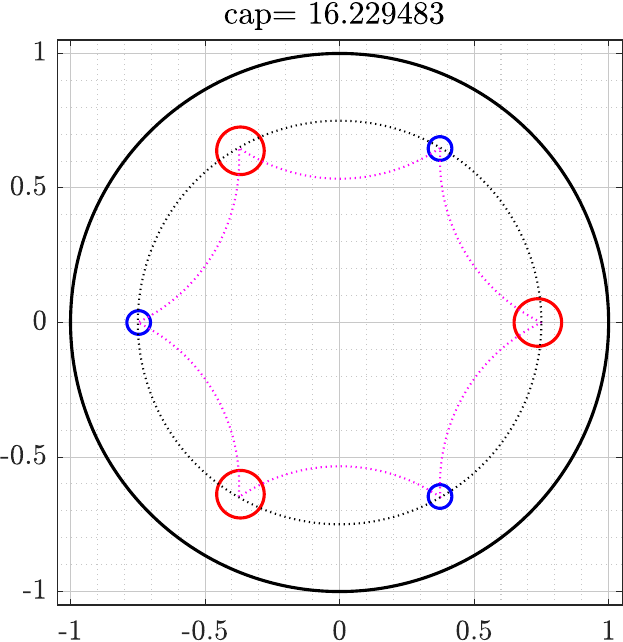}
	}\\
	{
		\includegraphics[width=0.31\textwidth]{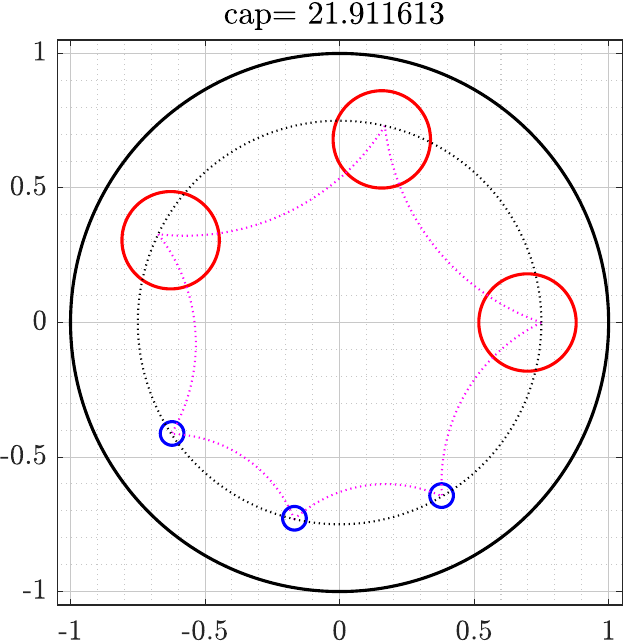}
	}
	{
		\includegraphics[width=0.31\textwidth]{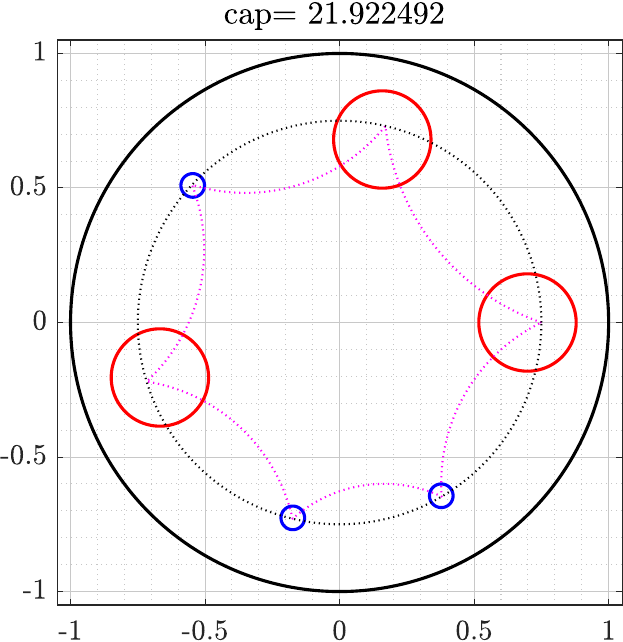}
	}
	{
		\includegraphics[width=0.31\textwidth]{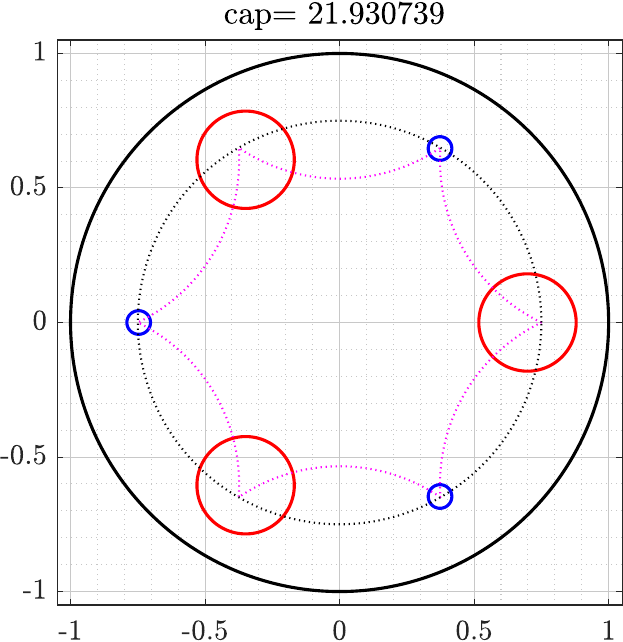}
	}	
	\caption{First row: Three disks with hyperbolic radii 0.2 and three disks with hyperbolic radii 0.4. Second row: Three disks with hyperbolic radii 0.2 and three disks with hyperbolic radii 0.8.}\label{fig:6Dn3}
\end{figure}

If the hyperbolic radii of two of these six disks are changed to either
$0.4$ or $0.8$, the natural symmetries induce three local maxima as shown in 
Figure~\ref{fig:6Dn2}. The hyperbolic distances between the centers of
any two adjacent disks for all cases in Figure~\ref{fig:6Dn2} are shown
in Table~\ref{tab:6Dn2}.
Similarly, with three disks three local maxima are observed 
(see Figure~\ref{fig:6Dn3} and Table~\ref{tab:6Dn3}).

Considering the results for constellations of disks with unequal
radii we can observe that in all cases  the maximal
dispersion property is again observed:
In the configuration which maximizes the capacity $\capa(\mathbb{B}^2,E)$ 
the positions of these six disks are on the Euclidean circle $|z|=R$.

\begin{table}
	\caption{The hyperbolic distances between the centers of consecutive disks in Figure~\ref{fig:6Dn1}.}
	\label{tab:6Dn1}
	\centering
	\begin{tabular}{lccccccc} 
		\hline
		Case & Capacity & $\rho(z_1,z_2)$ & $\rho(z_2,z_3)$ & $\rho(z_3,z_4)$ & $\rho(z_4,z_5)$ & $\rho(z_5,z_6)$ & $\rho(z_6,z_1)$ \\ \hline
		A & $13.7574$ & $2.6161$ & $2.6161$ & $2.6161$  & $2.6161$ & $2.6161$ & $2.6161$\\
		B & $14.6023$ & $2.7393$ & $2.5532$ & $2.5482$  & $2.5482$ & $2.5532$ & $2.7393$\\
		C & $16.6416$ & $2.9128$ & $2.4504$ & $2.4363$  & $2.4363$ & $2.4504$ & $2.9128$\\
		\hline
	\end{tabular}
\end{table}

\begin{figure}
	\centering
	\includegraphics[width=0.475\textwidth]{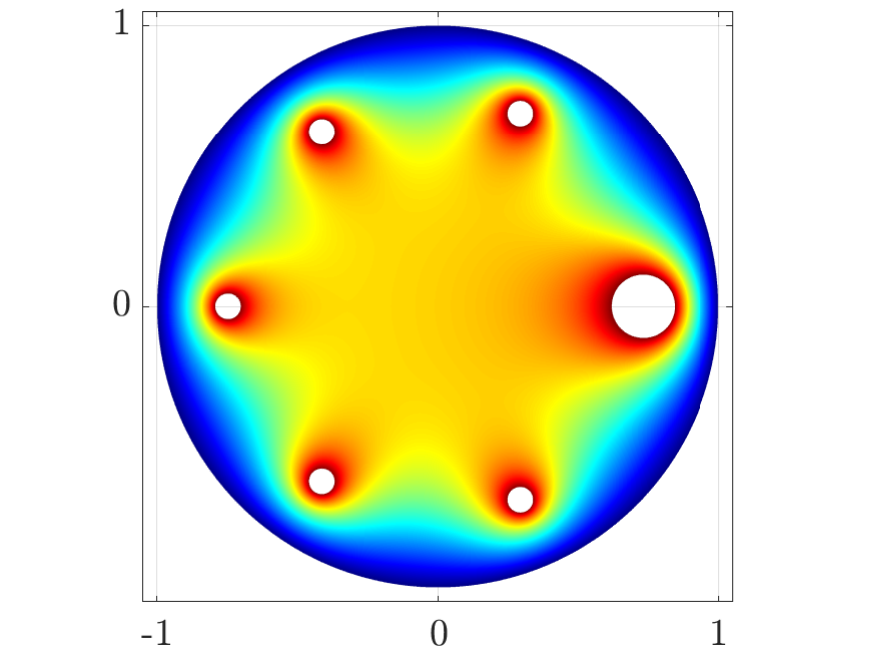} \hfill
	\includegraphics[width=0.475\textwidth]{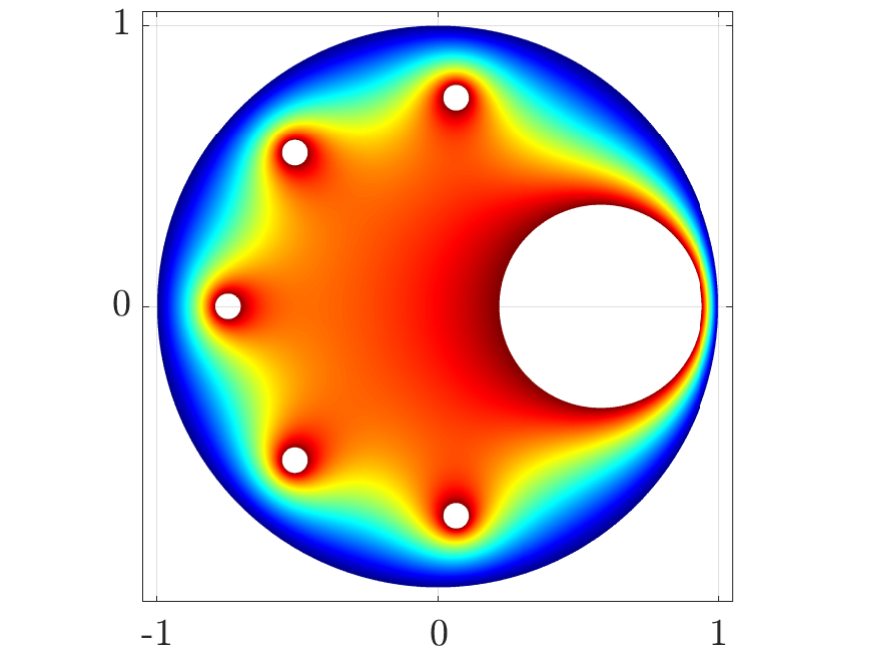} \\
	\includegraphics[width=0.475\textwidth]{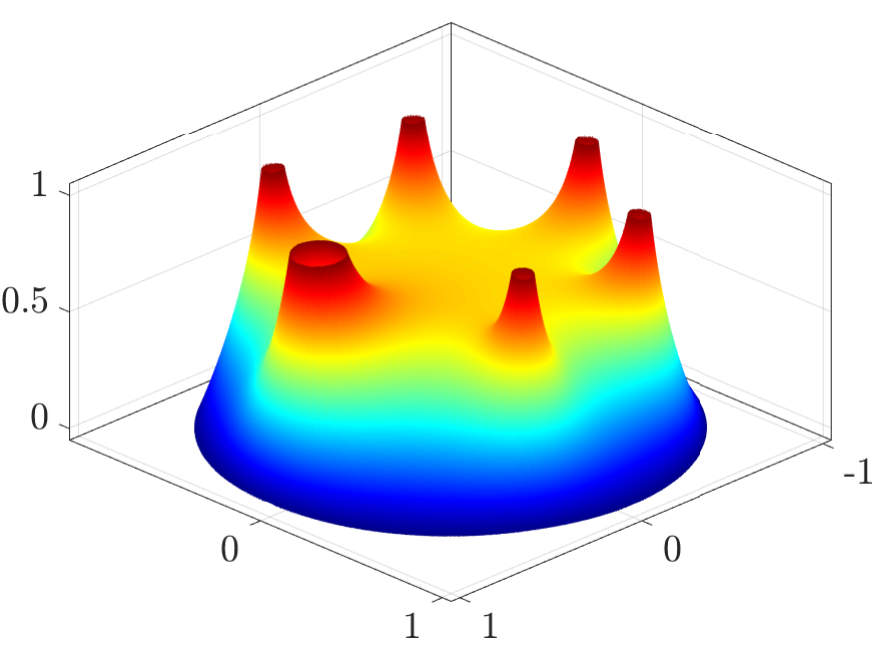} \hfill
	\includegraphics[width=0.475\textwidth]{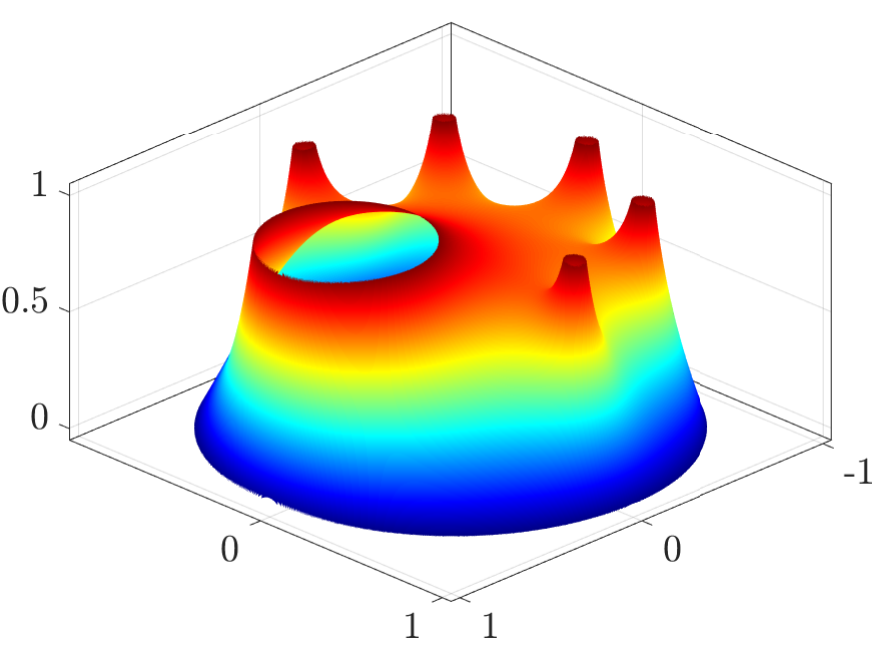} 
	\caption{Surface plots of the potentials for the six hyperbolic disks where the Disk 1 has hyperbolic radius $r_1$ and Disks 2 through 6 have hyperbolic radii $0.2$ for $r_1=0.5$ (left) and $r_1=1.5$ (right). The centers are at the positions that maximize the capacity.}\label{fig:6DPot}
\end{figure}

\begin{table}
	\caption{The hyperbolic distances between the centers of consecutive disks in Figure~\ref{fig:6Dn2}.}
	\label{tab:6Dn2}
	\centering
	\begin{tabular}{lccccccc} 
		\hline
		Case & Capacity & $\rho(z_1,z_2)$ & $\rho(z_2,z_3)$ & $\rho(z_3,z_4)$ & $\rho(z_4,z_5)$ & $\rho(z_5,z_6)$ & $\rho(z_6,z_1)$ \\ \hline
		A & $15.4245$ & $2.8523$ & $2.6767$ & $2.4861$  & $2.4815$ & $2.4861$ & $2.6767$\\
		B & $15.4263$ & $2.6800$ & $2.6800$ & $2.6742$  & $2.4867$ & $2.4867$ & $2.6742$\\
		C & $15.4266$ & $2.6747$ & $2.4920$ & $2.6747$  & $2.6747$ & $2.4920$ & $2.6747$\\
		D & $19.3499$ & $3.1452$ & $2.7534$ & $2.2784$  & $2.2672$ & $2.2784$ & $2.7534$\\
		E & $19.3620$ & $2.7700$ & $2.7700$ & $2.7519$  & $2.2813$ & $2.2813$ & $2.7519$\\
		F & $19.3643$ & $2.7553$ & $2.2965$ & $2.7553$  & $2.7553$ & $2.2965$ & $2.7553$\\
		\hline
	\end{tabular}
\end{table}

\begin{table}
	\caption{The hyperbolic distances between the centers of consecutive disks in Figure~\ref{fig:6Dn3}.}
	\label{tab:6Dn3}
	\centering
	\begin{tabular}{lccccccc} 
		\hline
		Case & Capacity & $\rho(z_1,z_2)$ & $\rho(z_2,z_3)$ & $\rho(z_3,z_4)$ & $\rho(z_4,z_5)$ & $\rho(z_5,z_6)$ & $\rho(z_6,z_1)$ \\ \hline
		A & $16.2261$ & $2.7931$ & $2.7931$ & $2.6126$  & $2.4210$ & $2.4210$ & $2.6126$\\
		B & $16.2280$ & $2.7903$ & $2.6189$ & $2.6161$  & $2.6108$ & $2.4263$ & $2.6135$\\
		C & $16.2295$ & $2.6161$ & $2.6161$ & $2.6161$  & $2.6161$ & $2.6161$ & $2.6161$\\
		D & $21.9116$ & $3.0047$ & $3.0047$ & $2.6025$  & $2.1205$ & $2.1205$ & $2.6025$\\
		E & $21.9225$ & $3.0015$ & $2.6190$ & $2.6172$  & $2.6035$ & $2.1336$ & $2.6053$\\
		F & $21.9307$ & $2.6161$ & $2.6161$ & $2.6161$  & $2.6161$ & $2.6161$ & $2.6161$\\
		\hline
	\end{tabular}
\end{table}

\subsection{Constellation of five disks constrained to the real line}
Next we consider a constellation of five hyperbolic disks under the constraint that 
the hyperbolic centers of these
disks lie within the interval $[-R,R]$.
The disks are numbered ($D_1,\ldots,D_5$) from left to right, and
$R=0.75$ in all experiments.

The set of experiments follows that of the previous section.
Four cases are considered: 
(a) all five disks have equal hyperbolic radii $= 0.2$,
(b) one of the disks has radius $= 0.4$, 
(c) two disks have radius $= 0.4$,
and finally 
(d) three disks have radius $= 0.4$.

\begin{figure}
	\centering
	{
		\includegraphics[width=0.31\textwidth]{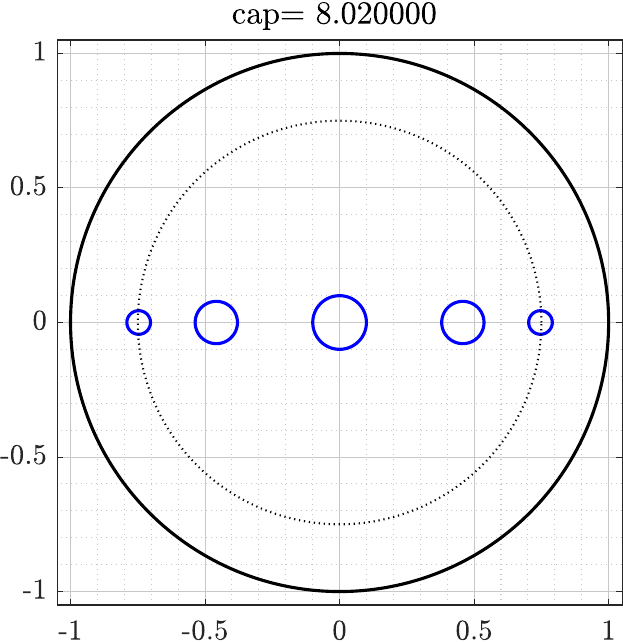}
	}\\
	{
		\includegraphics[width=0.31\textwidth]{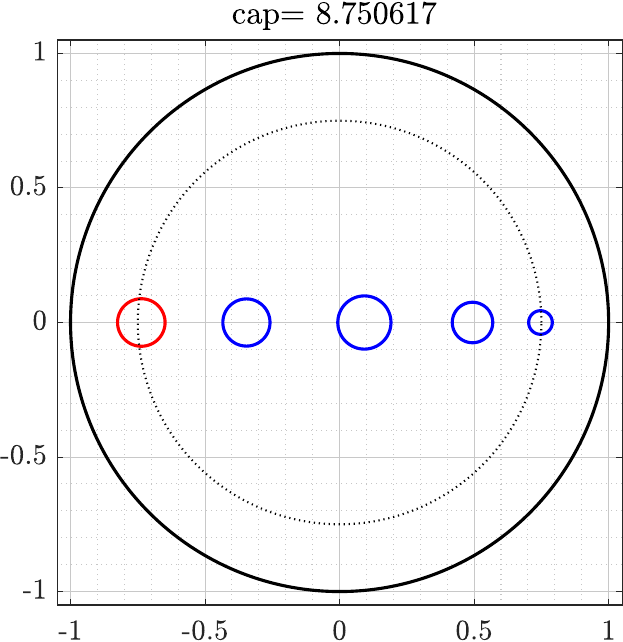}
	}
	{
		\includegraphics[width=0.31\textwidth]{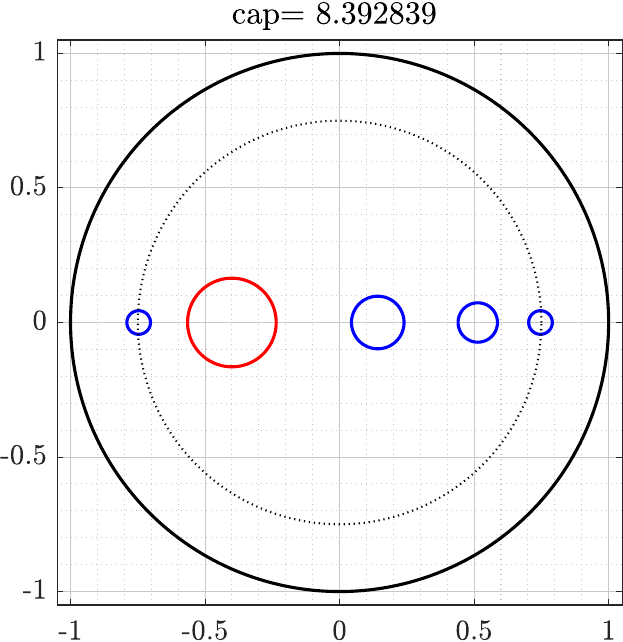}
	}
	{
		\includegraphics[width=0.31\textwidth]{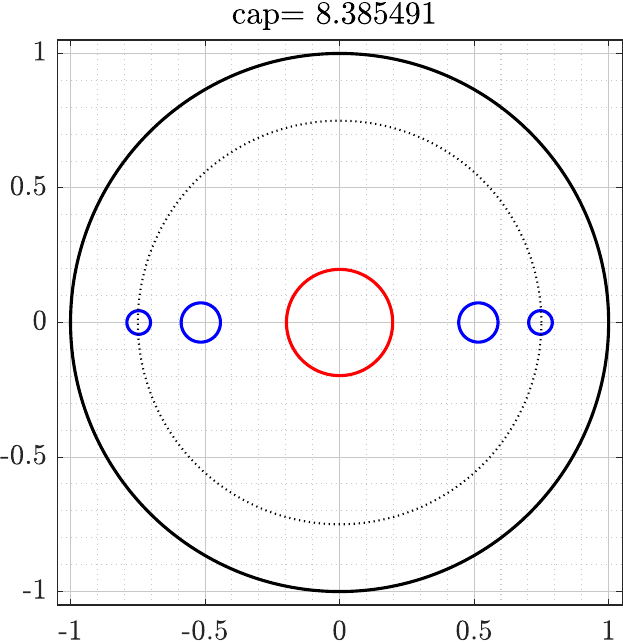}
	}
	\caption{First row: All disks have hyperbolic radius $0.2$. Second row: One disk with hyperbolic radius $0.4$ (red) and four disks with hyperbolic radii 0.2 (blue).}\label{fig:Re5Dn1}
\end{figure}

\begin{figure}
	\centering
	{
		\includegraphics[width=0.31\textwidth]{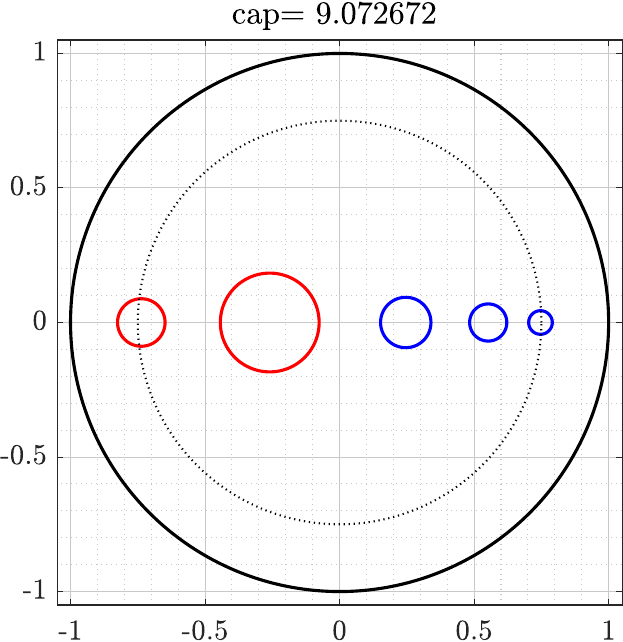}
	}
	{
		\includegraphics[width=0.31\textwidth]{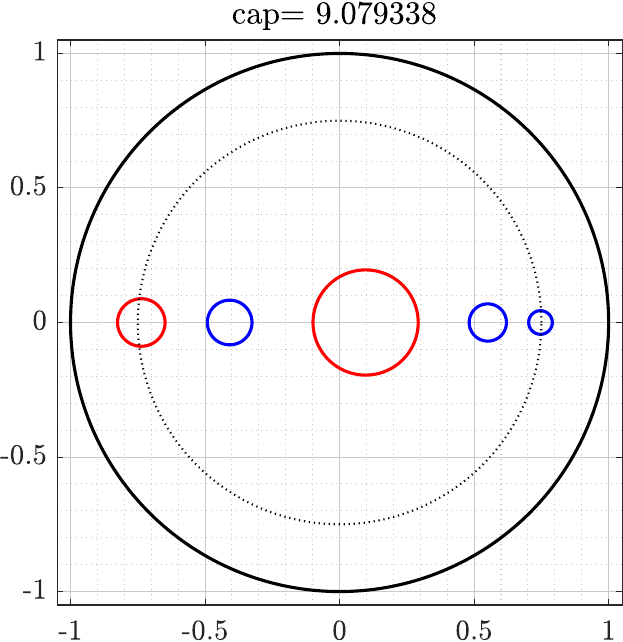}
	}
	{
		\includegraphics[width=0.31\textwidth]{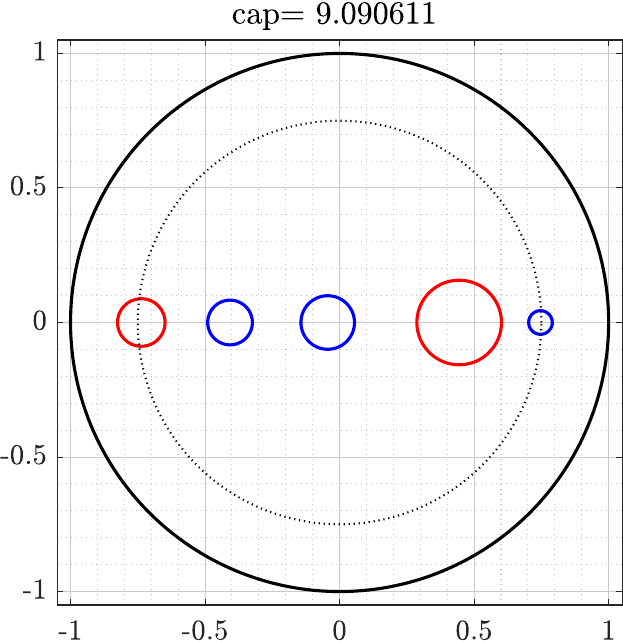}
	}\\
	{
		\includegraphics[width=0.31\textwidth]{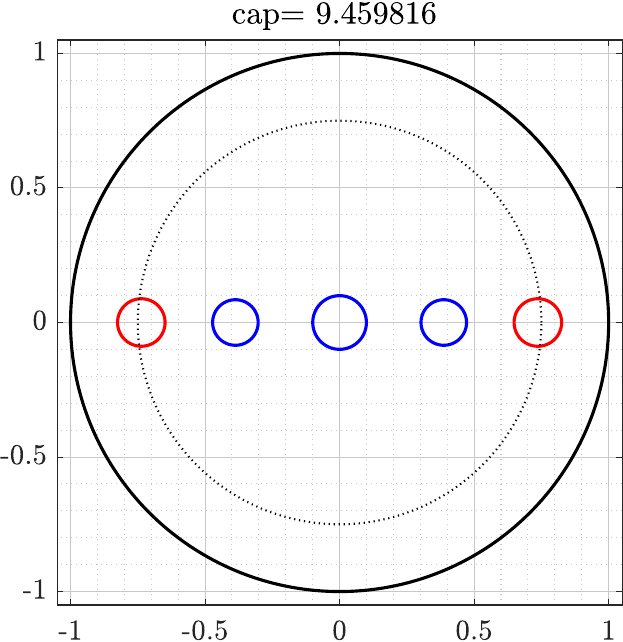}
	}
	{
		\includegraphics[width=0.31\textwidth]{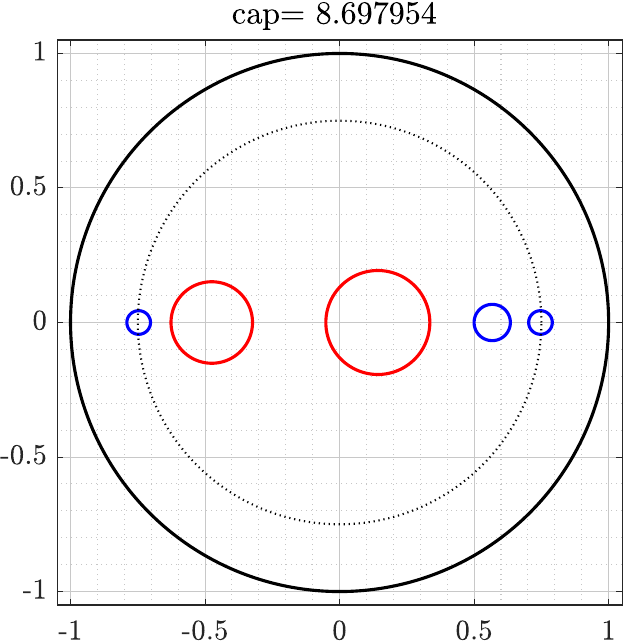}
	}
	{
		\includegraphics[width=0.31\textwidth]{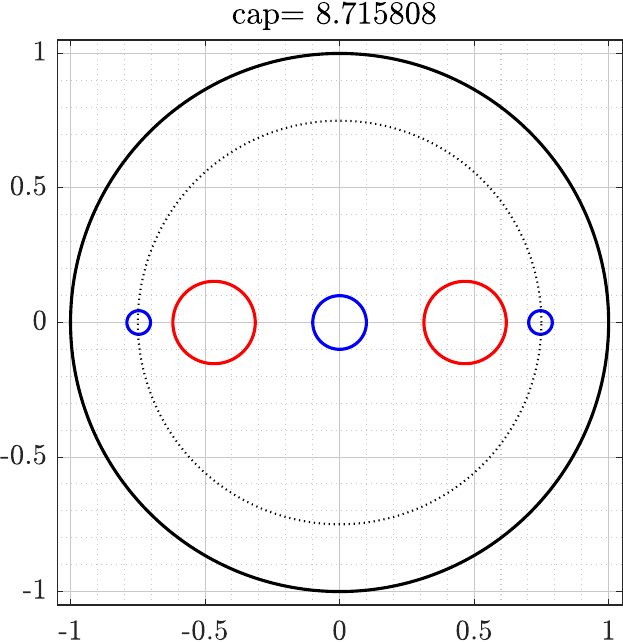}
	}
	\caption{Two disks have hyperbolic radius $0.4$ (red) and three disks have hyperbolic radii 0.2 (blue).}\label{fig:Re5Dn2}
\end{figure}

\begin{figure}
	\centering
	{
		\includegraphics[width=0.31\textwidth]{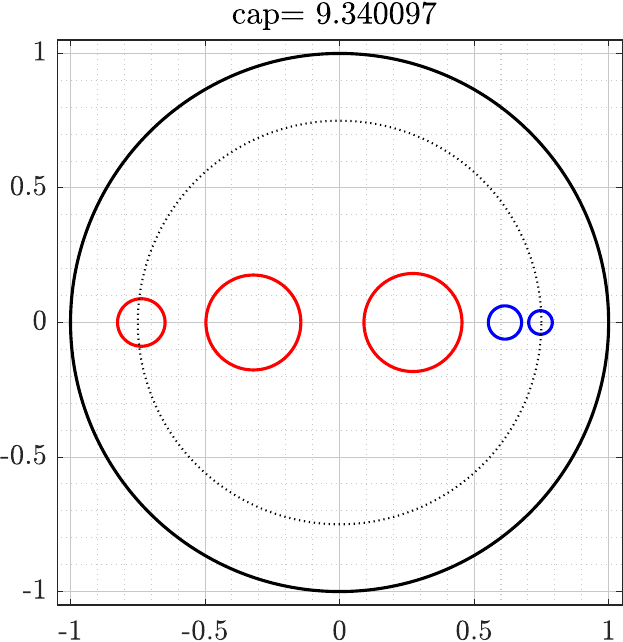}
	}
	{
		\includegraphics[width=0.31\textwidth]{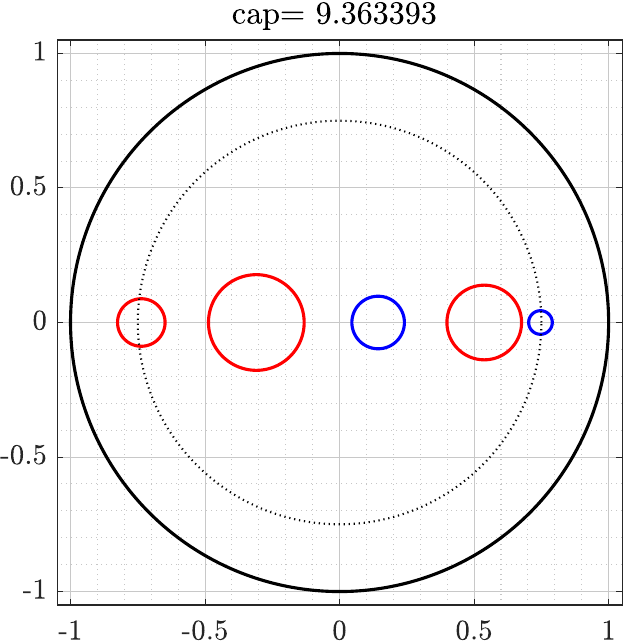}
	}
	{
		\includegraphics[width=0.31\textwidth]{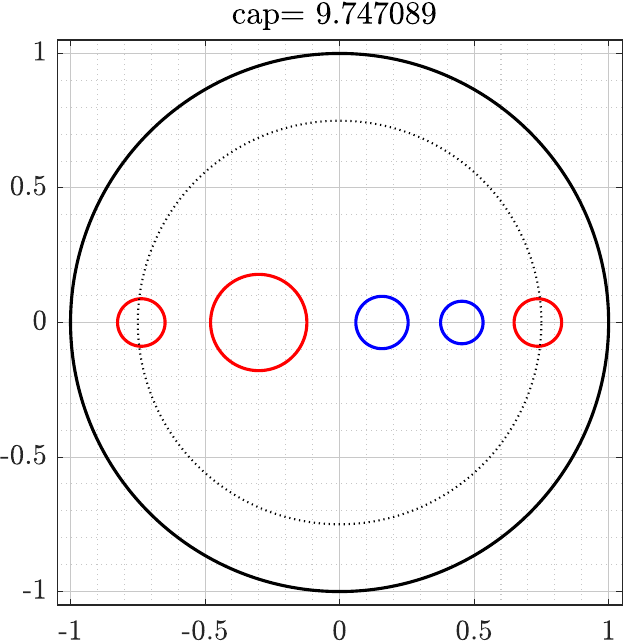}
	}\\
	{
		\includegraphics[width=0.31\textwidth]{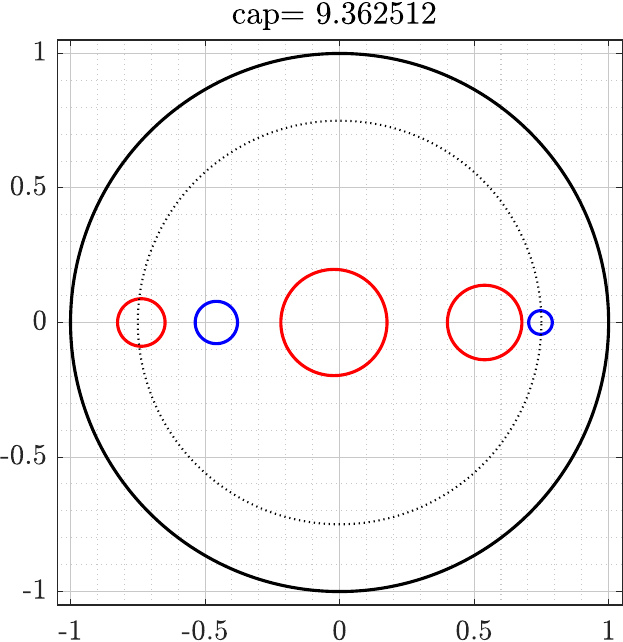}
	}
	{
		\includegraphics[width=0.31\textwidth]{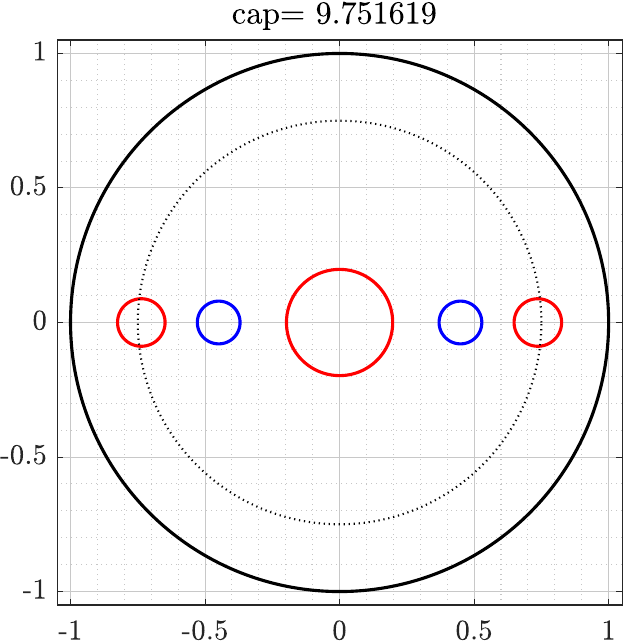}
	}
	\caption{Three disks have hyperbolic radius $0.4$ (red) and two disks have hyperbolic radii 0.2 (blue).}\label{fig:Re5Dn3}
\end{figure}

All configurations up to symmetry are summarized in 
Figures~\ref{fig:Re5Dn1}, ~\ref{fig:Re5Dn2}, ~\ref{fig:Re5Dn3},
and Tables~\ref{tab:Re5Dn1}, ~\ref{tab:Re5Dn2}, ~\ref{tab:Re5Dn3},
for (a) and (b), (c), and (d), respectively.
The maximal configurations exhibit the maximal dispersion property on a diameter:
$D_1$ and $D_5$ lie at the end points of the interval, $D_1$ and $D_5$ have the
largest radii, and if there are two or more disks with equal and largest radius, then
the distances between the disks are symmetric about the origin.

\begin{table}
	\caption{The hyperbolic distances between the centers of consecutive disks in Figure~\ref{fig:Re5Dn1}.}
	\label{tab:Re5Dn1}
	\centering
	\begin{tabular}{lccccc} 
		\hline
		Case & Capacity &$\rho(z_1,z_2)$ & $\rho(z_2,z_3)$  & $\rho(z_3,z_4)$ & $\rho(z_4,z_5)$  \\ \hline
		A & $8.0200$ & $0.9467$  & $0.9992$  & $0.99920$  & $0.9467$ \\
		B & $8.7506$ & $1.2166$  & $0.9160$  & $0.90604$  & $0.8532$ \\
		C & $8.3928$ & $1.0656$  & $1.1689$  & $0.85703$  & $0.8003$ \\
		D & $8.3855$ & $0.7943$  & $1.1516$  & $1.1516$   & $0.7943$ \\
		\hline
	\end{tabular}
\end{table}

\begin{table}
	\caption{The hyperbolic distances between the centers of consecutive disks in Figure~\ref{fig:Re5Dn2}.}
	\label{tab:Re5Dn2}
	\centering
	\begin{tabular}{lccccc} 
		\hline
		Case & Capacity & $\rho(z_1,z_2)$ & $\rho(z_2,z_3)$  & $\rho(z_3,z_4)$ & $\rho(z_4,z_5)$  \\ \hline
		A & $9.0727$ & $1.3935$  & $1.0596$ & $0.7477$  & $0.6910$ \\
		B & $9.0793$ & $1.0701$  & $1.0777$ & $1.0488$  & $0.6953$ \\
		C & $9.0906$ & $1.0733$  & $0.7834$ & $1.0806$  & $0.9545$ \\
		D & $9.4598$ & $1.1210$  & $0.8249$ & $0.8249$  & $1.1210$ \\
		E & $8.6980$ & $0.8745$  & $1.3693$ & $1.0007$  & $0.6473$ \\
		F & $8.7158$ & $0.8971$  & $1.0488$ & $1.0488$  & $0.8971$ \\
		\hline
	\end{tabular}
\end{table}

\begin{table}
	\caption{The hyperbolic distances between the centers of consecutive disks in Figure~\ref{fig:Re5Dn3}.}
	\label{tab:Re5Dn3}
	\centering
	\begin{tabular}{lccccc} 
		\hline
		Case & Capacity & $\rho(z_1,z_2)$ & $\rho(z_2,z_3)$  & $\rho(z_3,z_4)$ & $\rho(z_4,z_5)$  \\ \hline
		A & $9.3401$ & $1.2558$ & $1.2719$ & $0.7477$   & $0.5005$ \\
		B & $9.3634$ & $1.2810$ & $0.9565$ & $0.9534$   & $0.7010$ \\
		C & $9.7471$ & $1.3014$ & $0.9656$ & $0.6684$   & $0.9564$ \\
		D & $9.3625$ & $0.9473$ & $0.9550$ & $1.2930$   & $0.6966$ \\
		E & $9.7516$ & $0.9693$ & $0.9766$ & $0.9766$   & $0.9693$ \\
		\hline
	\end{tabular}
\end{table}

\subsection{Condensation of a constellation of $m$ disks  into one disk}\label{compare}
We study now the condensation of a constellation $E$ of $m$ hyperbolic disks with equal radii $r$ into the case of one hyperbolic disk constellation with equal capacity, and compare the hyperbolic area and perimeter of the original and the new constellation. That is, we assume that $E=\cup_{j=1}^{m}B_{\rho}(z_j,r)$ and we will find the value of $R$ such that $\capa(\BB,E)=\capa(\BB,B_{\rho}(0,R))$.
Recall first that the hyperbolic area and hyperbolic perimeter of a hyperbolic disk $B_{\rho}(z_j,r)$ are by \cite[Thm 7.2.2, p. 132]{b} 
\begin{equation}\label{hperim}
4\pi\sh^2\left(\frac{r}{2}\right) \quad {\rm and} \quad 2\pi\sh r, 
\end{equation} 
respectively. 

Let $c=\capa(\mathbb{B}^2,E)$, which will be approximated numerically using the above discussed BIE method. Since, by \eqref{cgqm_5.14}, 
\begin{equation}\label{diskcap} 
\capa(\mathbb{B}^2, B_{\rho}(0,R))= 2 \pi/\log(1/\th(R/2)) \,,
\end{equation}
the value of the radius $R$ of a single disk $B_{\rho}(0,R)$ with capacity equal to $c$ satisfies 
\[
2 \, \pi/\log(1/{\rm th}(R/2)) = c 
\]
and hence 
\begin{equation}\label{eq:md-R} 
R= 2\arth(e^{-2\pi/c})=-\log\th(\pi/c).
\end{equation}

\begin{figure}
\centering{\includegraphics[width=0.33\textwidth]{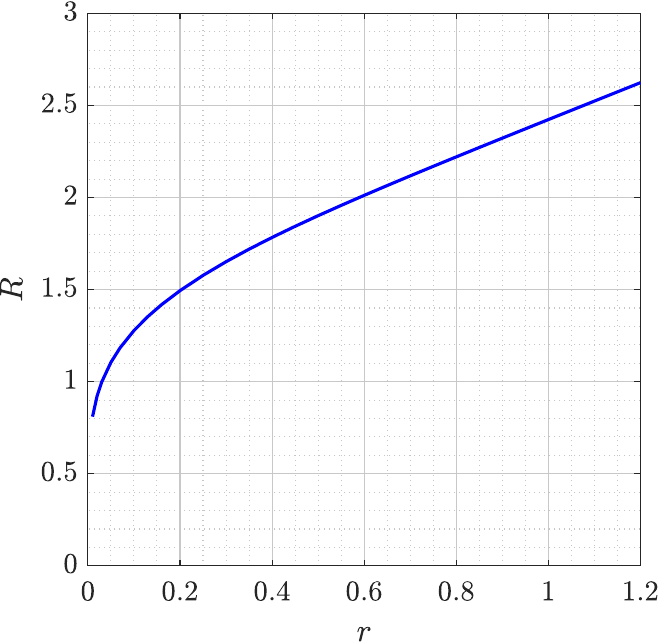}\hfill\includegraphics[width=0.33\textwidth]{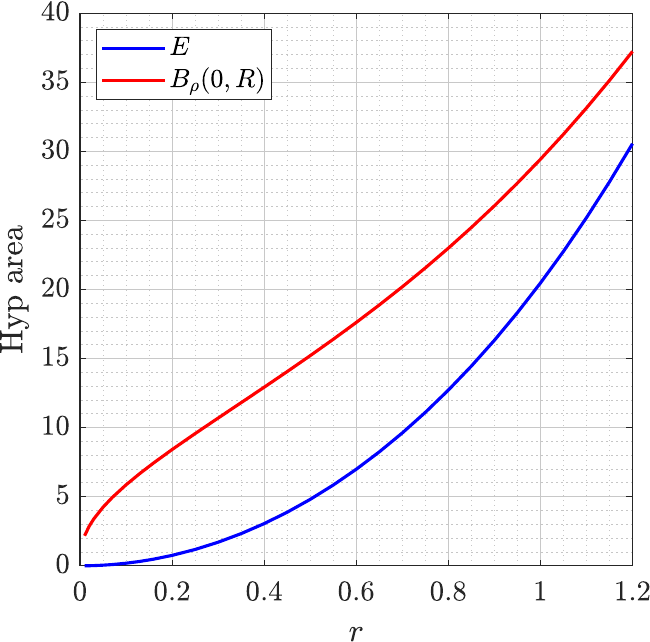}\hfill\includegraphics[width=0.33\textwidth]{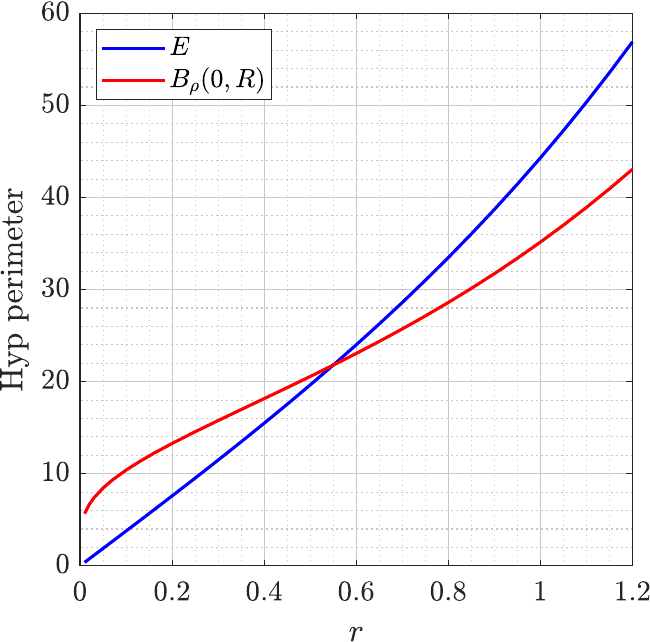}}
\caption{Left: The values of $R$ as a function of $r$ such that $\capa(\BB,B_{\rho}(0,R))=\capa(\BB,E)$ where $E=\cup_{j=1}^{m}B_{\rho}(z_j,r)$ for $m=6$ and $z_j=0.75e^{2\pi(j-1)\i/m}$, $j=1,2,\ldots,m$. Center and right: The hyperbolic area and perimeter of $E$ and $B_{\rho}(0,R)$ as functions of $r$.}\label{fig:md-RAP}
\end{figure}

As an example, we assume that $m=6$ and the centers $z_j$ of the hyperbolic disks $B_{\rho}(z_j,r)$ are given by
\[
z_j=0.75e^{2\pi(j-1)\i/m},\quad j=1,2,\ldots,m.
\]
We compute the capacity $c=\capa(\mathbb{B}^2,E)$ using the above BIE method with $n=2^{10}$. Then, we compute the values of $R$ via~\eqref{eq:md-R}. The computed values of $R$ for $0.1\le r\le 1.2$ are presented in Figure~\ref{fig:md-RAP}. Then, by~\eqref{hperim}, the hyperbolic area and perimeter of the disk $B_{\rho}(0,R)$ are equal to $4\pi\sh^2\left(\frac{R}{2}\right)$ and $2\pi\sh R$, respectively. Note that the hyperbolic area and perimeter of $E$ are given by
\[
4m\pi\sh^2\left(\frac{r}{2}\right) \quad {\rm and} \quad 2m\pi\sh r,
\]
respectively. The hyperbolic area and perimeter of $E$ and $B_{\rho}(0,R)$ are presented in Figure~\ref{fig:md-RAP}. The obtained results show that the hyperbolic area of the single disk $B_{\rho}(0,R)$ is always greater than the sum of the hyperbolic area of the six disks. However, the hyperbolic perimeter of the single disk $B_{\rho}(0,R)$ is greater than the sum of the hyperbolic perimeter of the six disks for small values of $r$. For large values of $r$, the perimeter of the six disks is greater than the perimeter of the single disk.

\section{Numerical Experiments: slit constellations}

In this section the elements of the constellations are hyperbolic segments of constant length.
The experiments follow the same pattern as those above, however, the constraints on configurations are
more restrictive. Again, we start with two segments and then increase complexity by adding more segments to the constellations.

\subsection{Constellation of two hyperbolic segments}

\begin{figure}
	\centering{
	\hfill\includegraphics[width=0.33\textwidth]{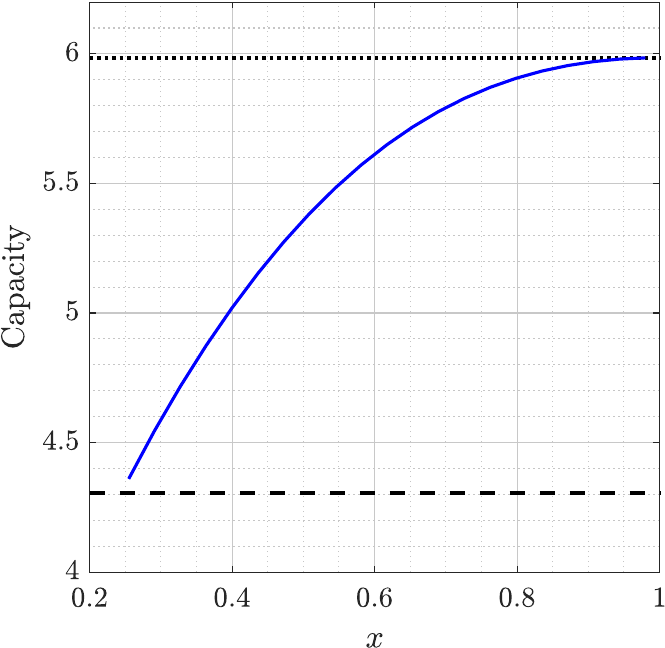} \hfill
	\includegraphics[width=0.33\textwidth]{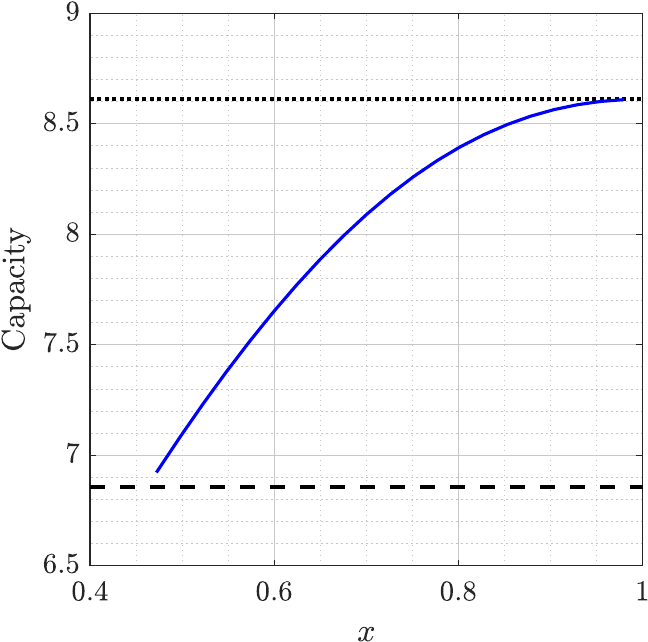}\hfill}
	\caption{Two hyperbolic segments with hyperbolic length $\ell=1$ (left) $\ell=2$ (right). The centers of these segments are $\pm x e^{\i\theta}$ for $\theta=0$.}\label{fig:S2Cap}
\end{figure}

We assume that the constellation $E$ is the union of two non-overlapping hyperbolic symmetric collinear segments $I_1$ an $I_2$ with equal hyperbolic length $\ell$ such that the centers of these segments are $\pm x e^{\i\theta}$ on the line $\arg(z)=\theta\in[0,\pi)$ where 
\[
\th\frac{\ell}{4}<x<1,
\] 
and hence $\rho\left(-xe^{\i\theta},xe^{\i\theta}\right)>\ell$.  The values of $\capa(\mathbb{B}^2,E)$ vs. $x$ are shown in Figure~\ref{fig:S2Cap} (left) for $\theta=0$, $\ell=1$ and in Figure~\ref{fig:S2Cap} (right) for $\theta=0$, $\ell=2$. 
Note that 
\[
\capa(\mathbb{B}^2,I_i)=\frac{2\pi}{\mu(\th(\ell/2))}, \quad i=1,2, 
\]
and hence 
\[
\capa(\mathbb{B}^2,I_1)+\capa(\mathbb{B}^2,I_2)=2\capa(\mathbb{B}^2,I_1)=\frac{4\pi}{\mu(\th(\ell/2))}
\]
is an upper bound for $\capa(\mathbb{B}^2,E)$. The values of this upper bound are shown in Figure~\ref{fig:S2Cap} as ``dotted line.''

The two segments merge into one segment $\hat I$ of hyperbolic length $2\ell$ when $x=\th\frac{\ell}{4}$. Thus
\[
\capa(\mathbb{B}^2,\hat I)=\frac{2\pi}{\mu(\th(\ell))} 
\]
is a lower bound for for $\capa(\mathbb{B}^2,E)$. The values of $\capa(\mathbb{B}^2,\hat I)$ are shown in Figure~\ref{fig:S2Cap} as ``dashed line.'' 

Figure~\ref{fig:S2Cap} shows that $\capa(\mathbb{B}^2,E)\to4\pi/\mu(\th(\ell/2))$ as $x\to1$ and $\capa(\mathbb{B}^2,E)\to2\pi/\mu(\th(\ell))$ as $x\to\th(\ell/4)$.

\subsection{Constellation of five radial hyperbolic segments with constant angle of separation}

\begin{figure}
	\centering{\hfill
	\includegraphics[width=0.33\textwidth]{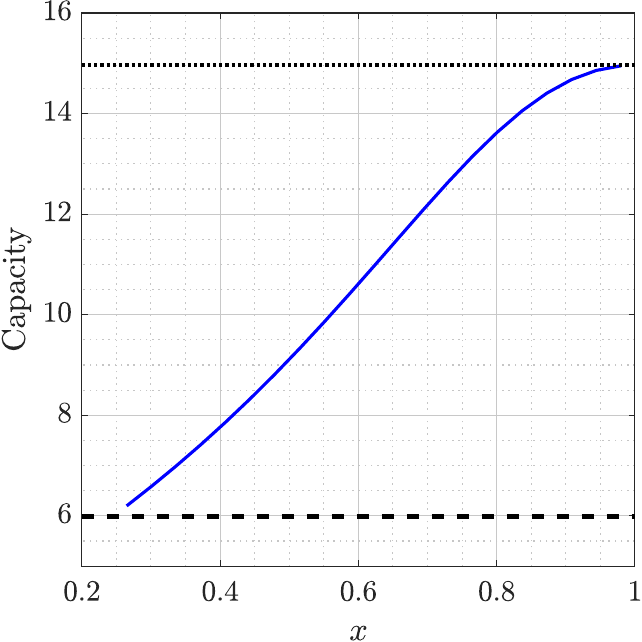} \hfill
	\includegraphics[width=0.33\textwidth]{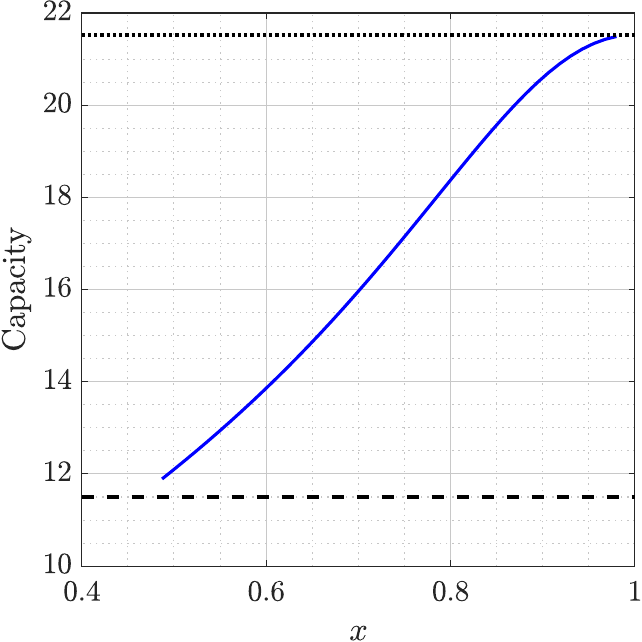}\hfill}
	\caption{Five radial hyperbolic segments with hyperbolic length $\ell=1$ (left) $\ell=2$ (right). The centers of these segments are $x e^{2k\pi\i/5}$, $k=1,\ldots,5$.}\label{fig:S5Cap}
\end{figure}

Next we let $E$ be the union of five non-overlapping hyperbolic segments, $I_1,\ldots,I_5$, with equal hyperbolic length $\ell$ such that the center of the segment $I_k$ is $x e^{2k\pi\i/5}$ where 
\[
\th\frac{\ell}{4}<x<1.
\] 
The computed approximate values of $\capa(\mathbb{B}^2,E)$ vs. $x$ are shown in Figure~\ref{fig:S5Cap} (left) for $\ell=1$ and in Figure~\ref{fig:S5Cap} (right) for $\ell=2$. 
Note that the five segments merge into one connected set $\hat I$ when $x=\th\frac{\ell}{4}$. Thus, using the same approach used in~\cite[Lemma 6.8]{nv-iso}, we can prove that
\[
\capa(\mathbb{B}^2,\hat I)=\frac{10\pi}{\mu(\th^5(\ell/2))} 
\]
which is a lower bound for $\capa(\mathbb{B}^2,E)$. The values of $\capa(\mathbb{B}^2,\hat I)$ are shown in Figure~\ref{fig:S5Cap} as ``dashed line.'' As in the previous example,  
\[
\sum_{k=1}^5\capa(\mathbb{B}^2,I_k)=5\capa(\mathbb{B}^2,I_1)=\frac{10\pi}{\mu(\th(\ell/2))}
\]
is an upper bound for $\capa(\mathbb{B}^2,E)$. The values of this upper bound are shown in Figure~\ref{fig:S5Cap} as ``dotted line.''

\begin{figure}
	\centering{\hfill
	\includegraphics[width=0.425\textwidth]{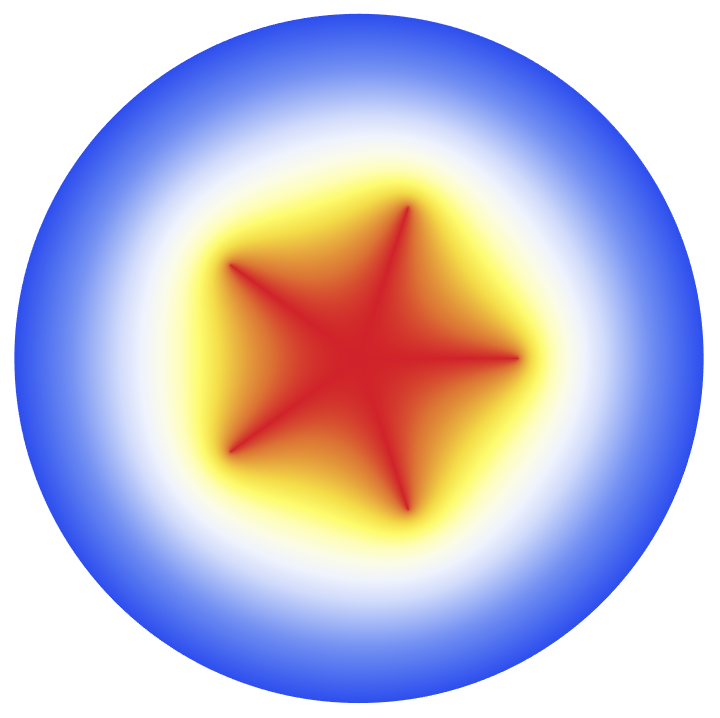} \hfill
	\includegraphics[width=0.425\textwidth]{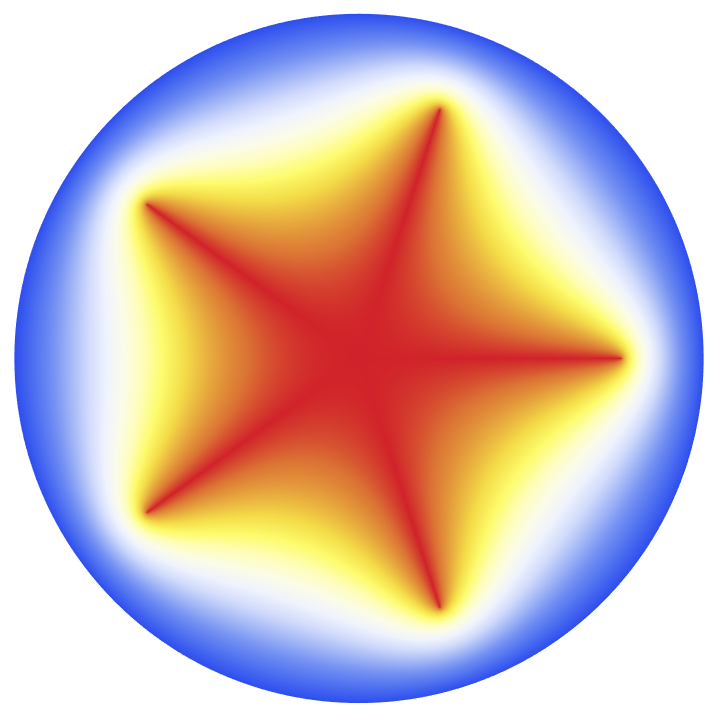}\hfill} \\
	\centering{\hfill\includegraphics[width=0.475\textwidth]{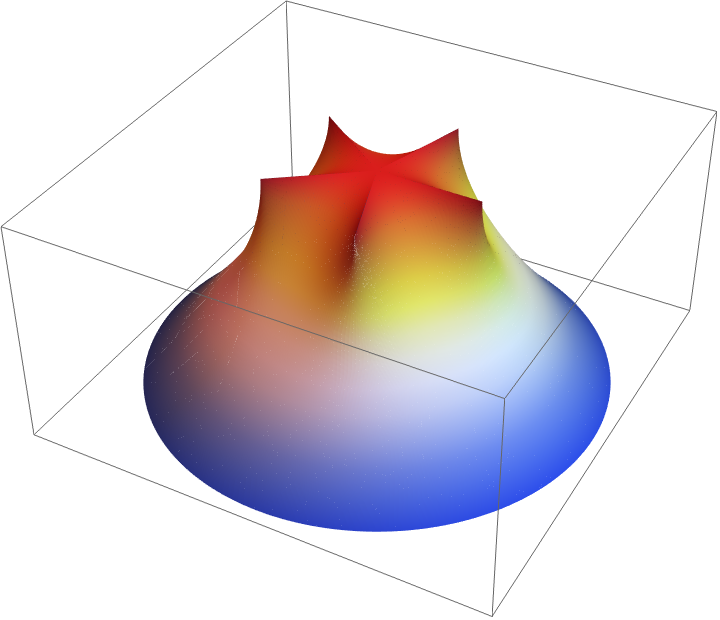} \hfill
	\includegraphics[width=0.475\textwidth]{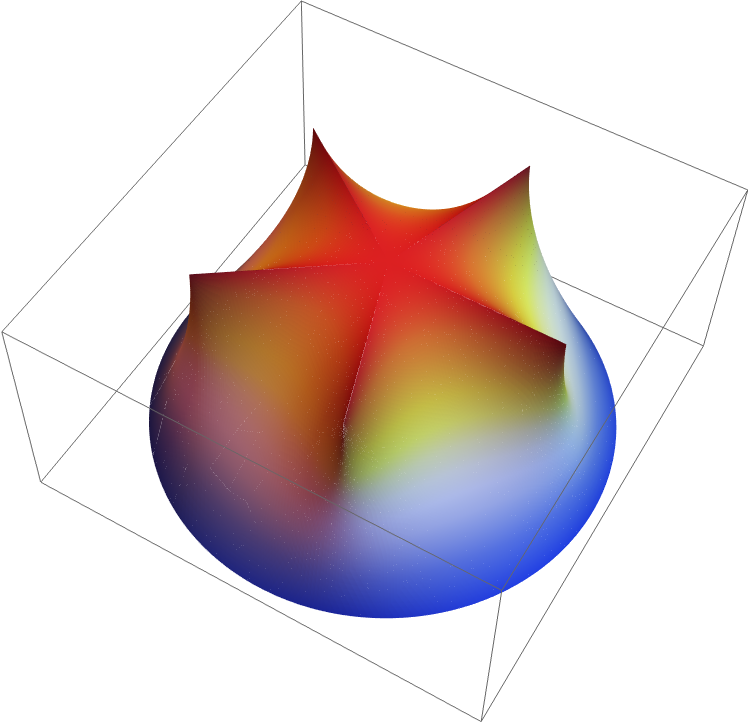}\hfill}
	\caption{Five radial hyperbolic segments. Surface plots of the potentials in the lower limit cases when
	all segments meet at the origin.}\label{fig:S5Pot}
\end{figure}

For numerical computing of the capacity $\capa(\mathbb{B}^2,\hat I)$, we use the $hp$-FEM where the absolute error in the computed capacity are $2\times 10^{-12}$ and $8\times 10^{-12}$ for the short and long segments, respectively. Plots of the potential function for the capacity $\capa(\mathbb{B}^2,E)$ are presented in Figure~\ref{fig:S5Pot}.

\begin{figure}
 \centering
 \subfloat{\includegraphics[height=2in]{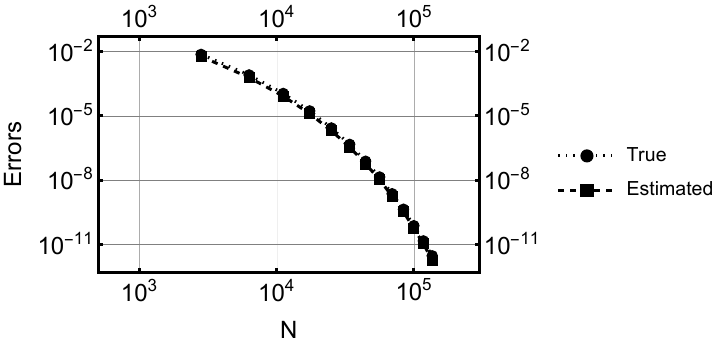}}\quad
 \subfloat{\includegraphics[height=2in]{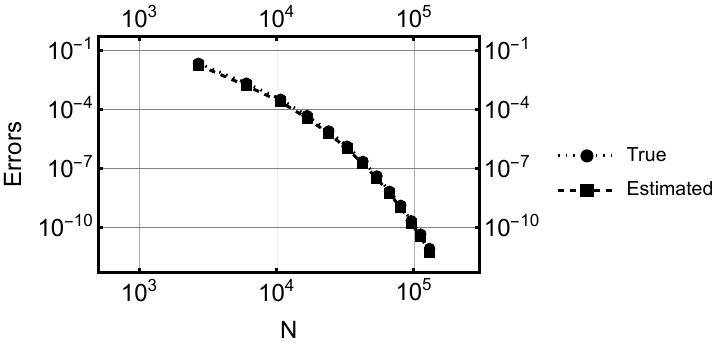}}
 \caption{Five radial hyperbolic segments. Error convergence: (Top) Short segments, (Bottom) Long segments.
 Both exact and estimated errors are shown in loglog-plots, $N$ is the number of degrees of freedom. }\label{fig:starconvergence}
\end{figure}

\subsection{Constellation of six hyperbolic segments constrained to a disk}

Analogously to the case with disks,
we consider the positions of six hyperbolic segments that maximize the capacity $\capa(\mathbb{B}^2,E)$ under the constraint that the hyperbolic centers of these disks are in the Euclidean disk $|z|\le R$ (we assume in the examples below that $R=0.75$). The segments are numbered $E_1$ to $E_6$ in counterclockwise orientation. We denote the center of the disk $E_j$ by $z_j$, $j=1,\ldots,6$.  Without loss of generality, we assume that the center $z_1$ of the segment $E_1$ is on the positive real axis. 

\begin{figure}
	\centering
	{
		\includegraphics[width=0.31\textwidth]{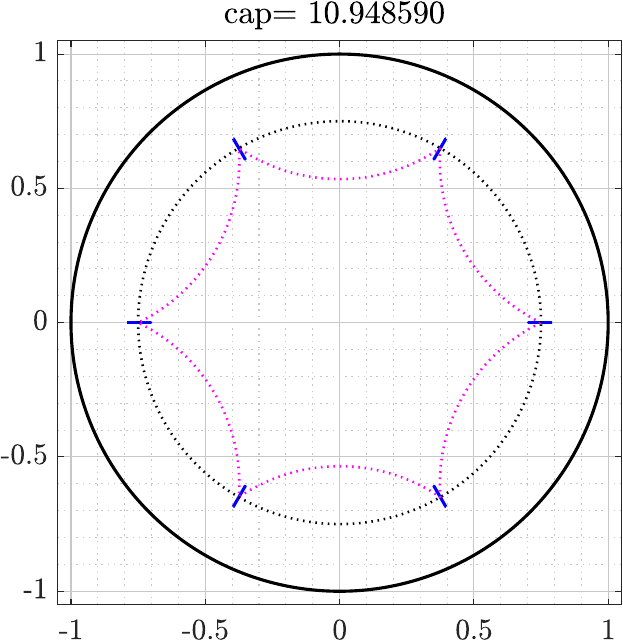}
	}
	{
		\includegraphics[width=0.31\textwidth]{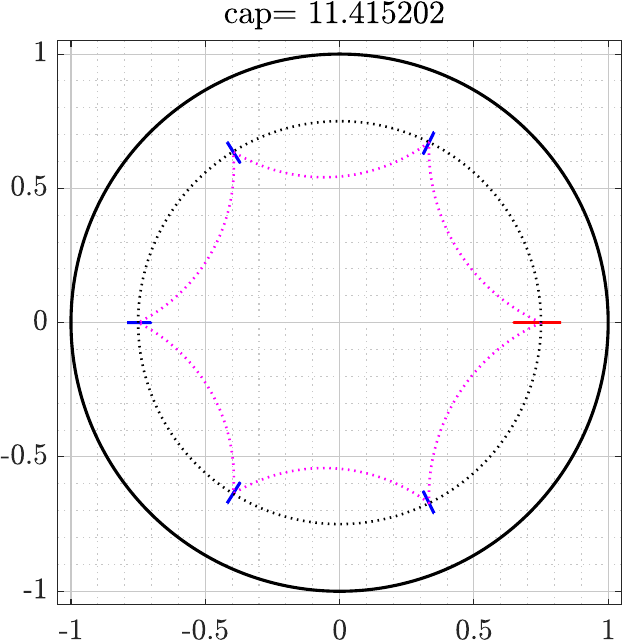}
	}
	{
		\includegraphics[width=0.31\textwidth]{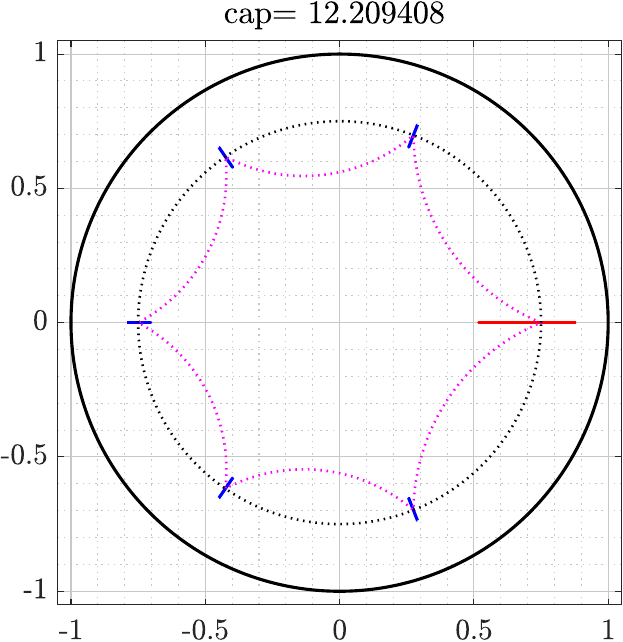}
	}
	\caption{Segment $E_1$ has hyperbolic length $0.4$ (left), $0.8$ (center), and $1.6$ (right). Segments $E_2$ through $E_6$ have hyperbolic length $0.4$.}\label{fig:6Sn1}
\end{figure}

First we assume that all six segments have equal hyperbolic length $= 0.4$. The positions of these six segments that maximize the capacity $\capa(\mathbb{B}^2,E)$ are on the Euclidean circle $|z|=R$ and such the hyperbolic distances between the centers of any two adjacent segments are equal (see Figure~\ref{fig:6Sn1} (left) and Table~\ref{tab:6Sn1}). When we change the hyperbolic length of one of these segments to be $0.8$ (see Figure~\ref{fig:6Sn1} (center)) or $1.6$ (see Figure~\ref{fig:6Sn1} (right)), then the centers of the other segments are moved away from the larger segment (see Table~\ref{tab:6Sn1}).

\begin{table}
	\caption{The hyperbolic distance between the centers of consecutive segments in Figure~\ref{fig:6Sn1}.}
	\label{tab:6Sn1}
	\centering
	\begin{tabular}{lccccccc} 
		\hline
		Case & Capacity & $\rho(z_1,z_2)$ & $\rho(z_2,z_3)$ & $\rho(z_3,z_4)$ & $\rho(z_4,z_5)$ & $\rho(z_5,z_6)$ & $\rho(z_6,z_1)$ \\ \hline
		A & $10.9486$ & $2.6161$ & $2.6161$ & $2.6161$  & $2.6161$ & $2.6161$ & $2.6161$\\
		B & $11.4152$ & $2.7080$ & $2.5705$ & $2.5655$  & $2.5655$ & $2.5705$ & $2.7080$\\
		C & $12.2094$ & $2.8236$ & $2.5089$ & $2.4931$  & $2.4931$ & $2.5089$ & $2.8236$\\
		\hline
	\end{tabular}
\end{table}

\subsection{Constellation of five hyperbolic segments constrained to the real line}

In the final experiment we consider the positions of five hyperbolic segments that maximize the capacity $\capa(\mathbb{B}^2,E)$ under the constraint that the hyperbolic centers of these slits are in the interval $[-R,R]$ (we assume in the examples below that $R=0.75$). The segments are numbered ($E_1,\ldots,E_5$) from left to right.

\begin{figure}
	\centering
	{
		\includegraphics[width=0.31\textwidth]{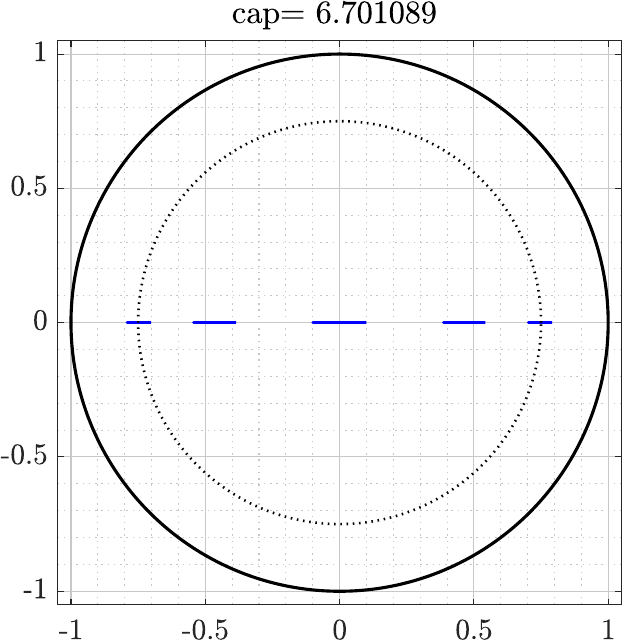}
	}\\
	{
		\includegraphics[width=0.31\textwidth]{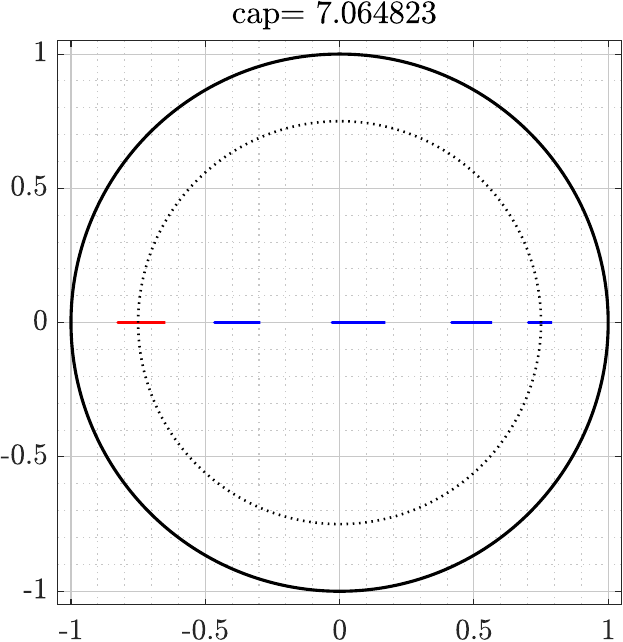}
	}
	{
		\includegraphics[width=0.31\textwidth]{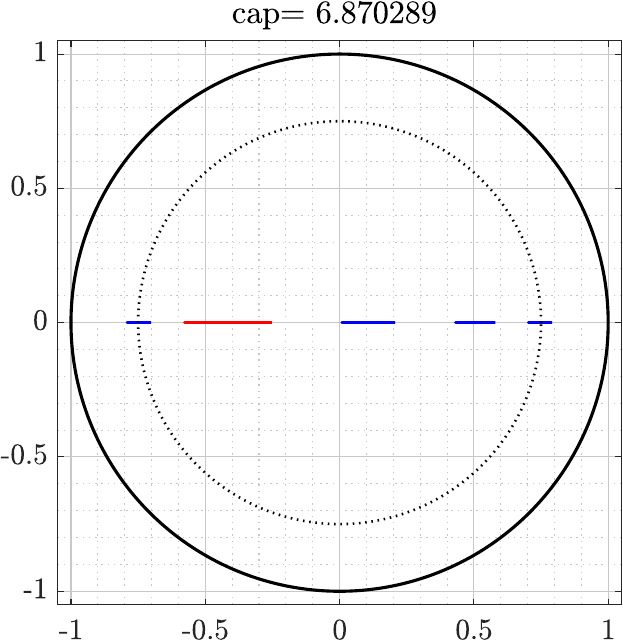}
	}
	{
		\includegraphics[width=0.31\textwidth]{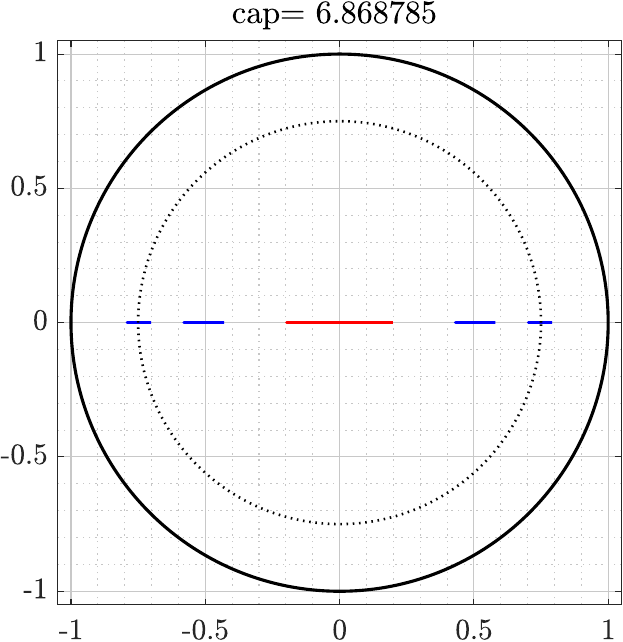}
	}
	\caption{First row: All segments have hyperbolic length $0.4$. Second row: One segment with hyperbolic length $0.8$ (red) and four segments with hyperbolic length $0.4$ (blue).}\label{fig:Re5Sn1}
\end{figure}

First all five segments are set to have equal hyperbolic length $= 0.4$. The positions of these five segments that maximize the capacity $\capa(\mathbb{B}^2,E)$ are shown in Figure~\ref{fig:Re5Sn1} and the hyperbolic distance between the centers of any two adjacent segments is presented in Table~\ref{tab:Re5Sn1}. Then we change the hyperbolic length of one of these segments to be $0.8$. The obtained results are presented in Figure~\ref{fig:Re5Sn1} and Table~\ref{tab:Re5Sn1}. 

\begin{table}
	\caption{The hyperbolic distances between the centers of consecutive segments in Figure~\ref{fig:Re5Sn1}.}
	\label{tab:Re5Sn1}
	\centering
	\begin{tabular}{lccccc} 
		\hline
		Case & Capacity & $\rho(z_1,z_2)$ & $\rho(z_2,z_3)$  & $\rho(z_3,z_4)$ & $\rho(z_4,z_5)$  \\ \hline
		A & $6.7011$ & $0.9293$  & $1.0166$ & $1.0166$  & $0.9293$ \\
		B & $7.0648$ & $1.1379$  & $0.9521$ & $0.9439$  & $0.8579$ \\
		C & $6.8703$ & $1.0317$  & $1.1313$ & $0.9071$  & $0.8218$ \\
		D & $6.8688$ & $0.8222$  & $1.1237$ & $1.1237$  & $0.8222$ \\
		\hline
	\end{tabular}
\end{table}

In all cases the results computed with BIE and FEM agree within the prescribed tolerance.

\subsubsection{On Computational Costs}
Optimization with the object function obtained by means of PDE solution is inherently
expensive.
In Table~\ref{tbl:timing} performance data on
the six disks maximization problem shown in Figure~\ref{fig:6Dint} is presented.
In all cases the interior-point tolerance is the same, $\epsilon = 10^{-6}$, and within the
$hp$-FEM simulations, meshing is performed with the same discretization
control in every evaluation. Not surprisingly, the overall conclusions are very similar
to those drawn in our previous work \cite{hnv}, where minimization was considered.
Comparison of the two methods is only qualitative, since both underlying
hardware and the interior-point implementations are different.

The two implementations have very different requirements per iteration step.
It is very likely that this is due to the fact that the numerical differentiation
algorithm of Matlab is different from the one of Mathematica.
Observe that the number of iteration steps becomes comparable once 
the $hp$-solutions are sufficiently accurate, yet the
number of evaluations is not. 
The average time for one evaluation in BIE
is four to five times faster than one evaluation in $hp$-FEM.
Matlab and Mathematica results have been computed on Lenovo ThinkPad with Intel Core i9-12900HX, 2300 Mhz, 16 Core(s), 24 Logical Processors and 64 GB RAM, and
Apple Silicon Mac mini M2 Pro 32GB (2023), respectively.
	
In short, for optimal performance, the individual solutions must be accurate
enough so that the error induced by numerical approximation
of the gradients and Hessians is balanced with other sources of error.
For BIE, the problem is practically fully resolved already at $n=2^6$,
whereas
for the $hp$-FEM it appears that the same mesh with $p=4$ is not adequate
in comparison with the one at $p=8$. Even though the time spent in one individual
iteration step is doubled, the overall time for $p=8$ is
significantly lower. 

\begin{table}[h]
\centering
\caption{Solution times for the maximization process for the six hyperbolic disks in Figure~\ref{fig:6Dint}. Number of steps is number of iterations in the interior-point algorithm. Number of evaluations is the total number of solves performed during the maximization.
	}\label{tbl:timing}
	\begin{tabular}{llrrr}
		Method & Discretization & Time & \# of steps & \# of evaluations \\\hline
		BIE      & $n=2^6$    & 162.5   & 16  & 204\\
		         & $n=2^8$    & 274.7   & 17  & 216 \\
		         & $n=2^{10}$ & 637.8   & 22  & 286 \\\hline
		$hp$-FEM & $p=4$      & 21749.2 & 144 & 23706 \\
		         & $p=6$      & 5132.6  & 24  & 3704  \\
		         & $p=8$      & 2170.9  & 6   & 1050   \\\hline
	\end{tabular}
\end{table}

\begin{figure}
	\centering
	{
		\includegraphics[width=0.48\textwidth]{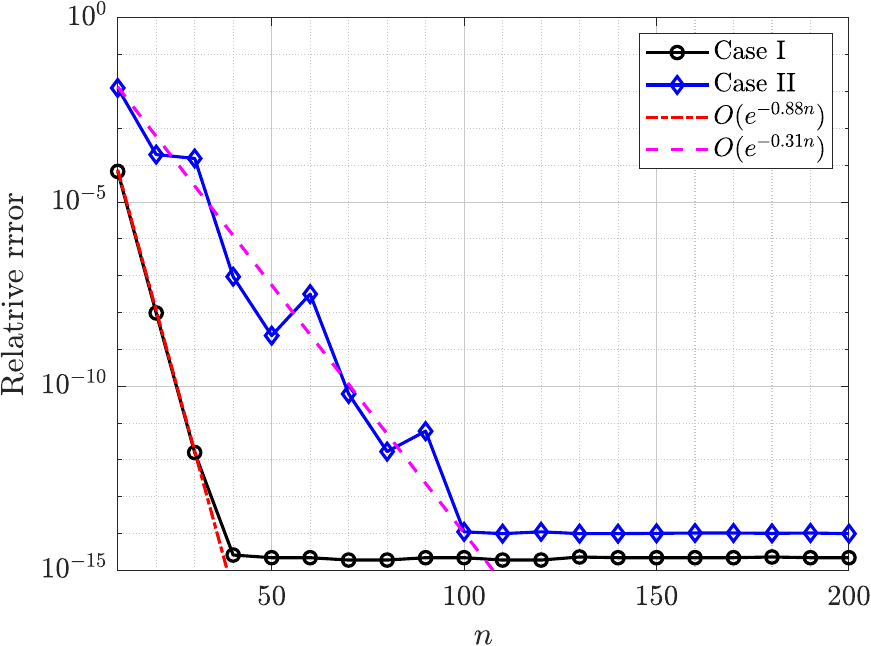}
	}
	{
		\includegraphics[width=0.48\textwidth]{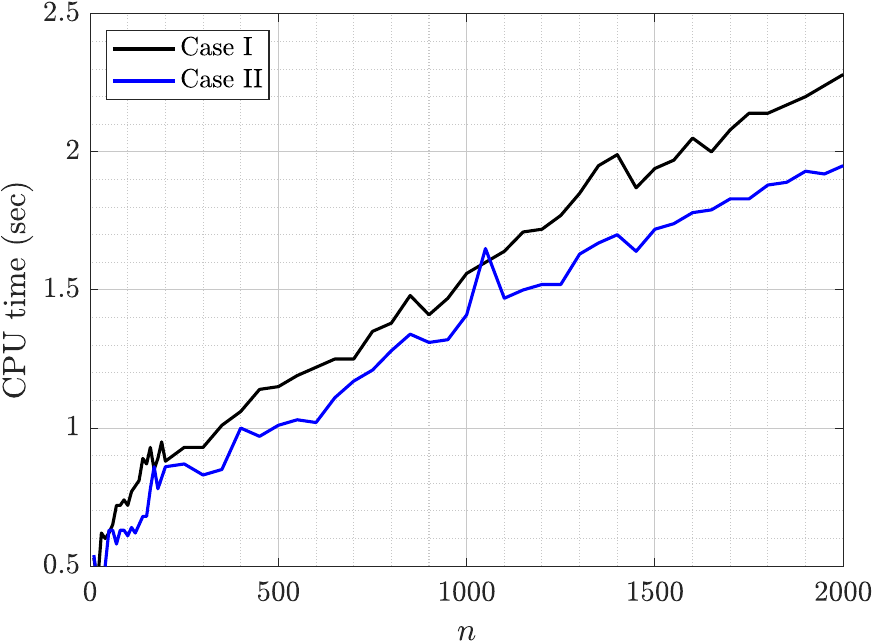}
	}
\caption{Left: The errors in the computed values of ${\rm cap}(\mathbb{B}^2,E)$ using the BIE method as functions of $n$ for both domains in Figure~\ref{fig:6Dint}.
Right: Run time in seconds as a function of $n$ required to compute ${\rm cap}(\mathbb{B}^2,E)$ using the BIE method.}
\label{fig:errBIE}
\end{figure}

\begin{figure}
	\centering
	{
		\includegraphics[height=1.5in]{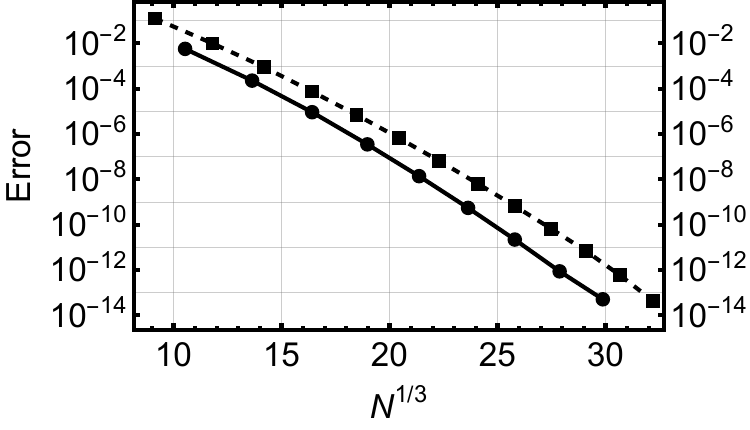}
	}
{
	\includegraphics[height=1.5in]{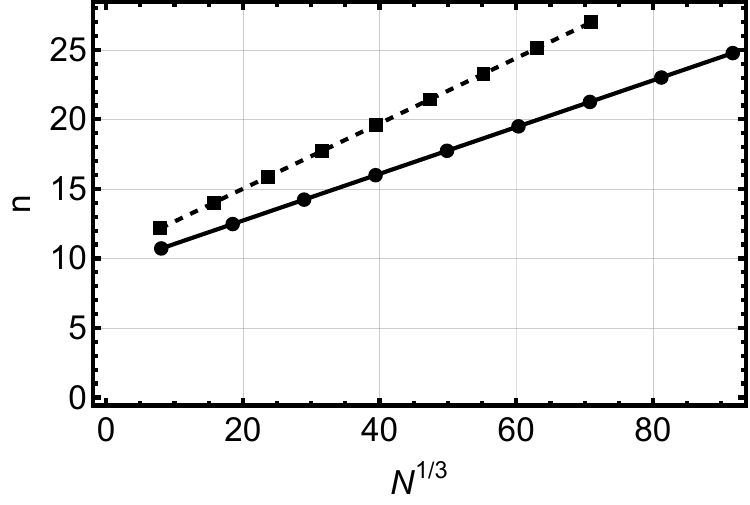}
}
\caption{%
{Left: The errors in the computed values of ${\rm cap}(\mathbb{B}^2,E)$
using the $hp$-FEM as functions of $N^{1/3}$ for both domains in
Figure~\ref{fig:6Dint}.
Right: Comparison of $N^{1/3}$ ($hp$-FEM) vs $n$ (BIE). The dashed line refers to the
maximal constellation.}}
\label{fig:errFEM}
\end{figure}

It is clear from Table~\ref{tbl:timing} that obtaining the positions of
the inner disks that maximize the capacity requires computing the
capacity for hundreds of times which can be computed accurately using
the two proposed methods. To illustrate the accuracy and the
computational cost of the two methods for each evaluation of the
capacity, we discuss in details the computation of the capacity of the
initial domain in Figure~\ref{fig:6Dint} (left) and the computed domain
with maximum capacity in Figure~\ref{fig:6Dint} (right).
For both domains, the value of the capacity is computed with $n=2^{13}$
to obtain ${\rm cap}(\mathbb{B}^2,E)=7.288440106260681$ for the domain in
Figure~\ref{fig:6Dint} (left) and ${\rm cap}(\mathbb{B}^2,E)=13.757382935965428$
for the domain in Figure~\ref{fig:6Dint} (right). 
These two computed values are considered as the reference values and used to estimate the error
in the values of ${\rm cap}(\mathbb{B}^2,E)$ computed by the two numerical methods. The computed error for the BIE method as a function of $n$ is presented in
Figure~\ref{fig:errBIE} (left) where Case I is the domain in
Figure~\ref{fig:6Dint} (left) and Case II is the domain in
Figure~\ref{fig:6Dint} (right). 
Similar convergence graphs are shown for the $hp$-FEM in Figure~\ref{fig:errFEM}.
For the cases with smooth boundaries both methods convergence exponentially
and we conclude that the complexity comparison can be reduced to
$n \sim N^{1/3}$, where $n$ is the number of quadrature points in BIE method
and $N$ is the number of degrees of freedom in FEM. Notice that this
rate for FEM is the optimal one given in Theorem~\ref{propermesh}.

\begin{figure}
	\centering
	{
		\includegraphics[width=0.48\textwidth]{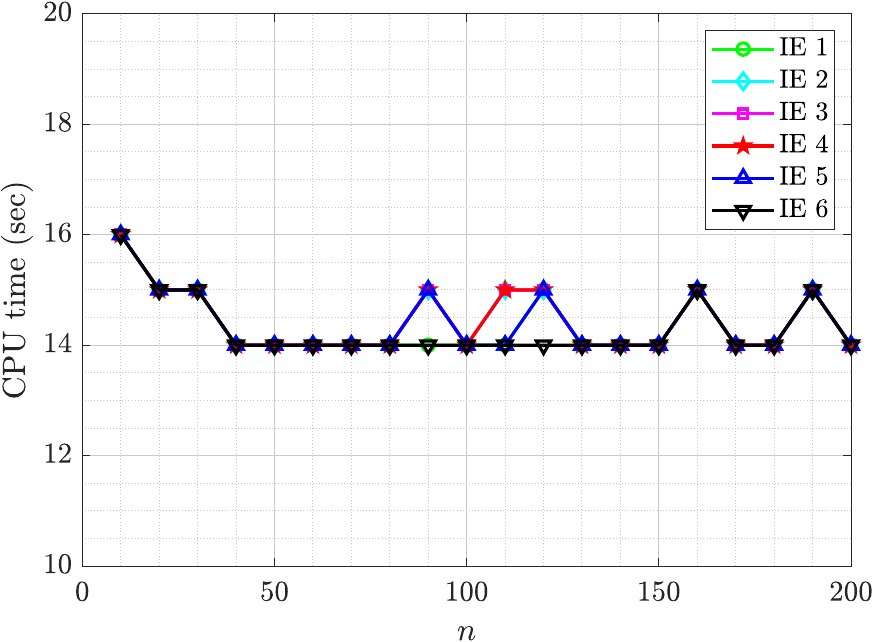}
	}
	{
		\includegraphics[width=0.48\textwidth]{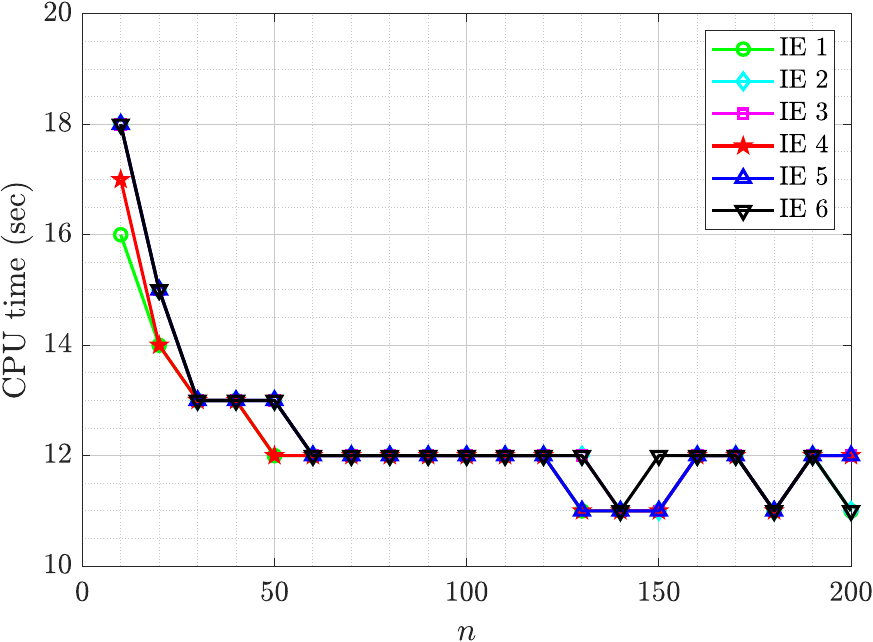}
	}
\caption{The number of GMRES iterations as a function of $n$ required to solve the $m=6$ integral equations for the two domains in Figure~\ref{fig:6Dint}. }
\label{fig:gmr}
\end{figure}

Computing the capacity ${\rm cap}(\mathbb{B}^2,E)$ using the BIE method requires solving the boundary integral equation~\eqref{eq:ie} for $m$ times. 
The computation of the right-hand side of the integral equation requires $O((m+1)n\ln n)$ operations and each iteration of the GMRES method requires $O((m+1)n)$ operations~\cite{Nas-ETNA}. Thus, solving the integral equation requires $O((m+1)n\ln n)$ operations and hence computing the capacity ${\rm cap}(\mathbb{B}^2,E)$ requires, approximately, $O((m+1)^2n\ln n)$ operations. 
To corroborate this, we present in Figure~\ref{fig:errBIE} (right) a plot of the CPU time (in seconds) as a function of $n$ required to compute the capacity of the two domains in Figure~\ref{fig:6Dint}. It is clear from these graphs that the CPU time depends almost linearly on $n$. We also provide in Figure~\ref{fig:gmr} the number of GMRES iterations required to solve the $m=6$ integral equations for the two domains. These graphs show that the number of GMRES iterations is almost independent of $n$. On one hand, the number of iterations for the domain in Figure~\ref{fig:6Dint} (left) is less than the number of iterations for the domain in Figure~\ref{fig:6Dint} (right) since the disks on the right are well-separated compared to the disks on the left. This implies that the CPU time required to compute the capacity of the domain on the right is less than the CPU time for the domain on the left. On the other hand, it is clear from Figure~\ref{fig:errBIE} that the order of convergence for the domain on the left is better than the order of convergence for the domain on the right. A possible reason for that is the inner disks for the domain on the right are closer to the external boundary compared to the domain on the left.
A detailed analysis of the order of convergence for the above BIE method is beyond the scope of the current paper and will be a subject of future work.

The computational complexity of the $hp$-FEM is comprehensively documented
in the references by Schwab \cite{schwab} and Szabo and Babu{\v{s}}ka \cite{sb}.
In contrast to the standard $h$-version the numerical integration of the inner products typically amounts to one half of
the solution times.

\begin{remark}
In comparison with a similar minimization problems in \cite{hnv}, we observe
that the constrained maximization problems are less resource intensive in terms of iteration steps
and runtimes. This is more notable in the BIE results. 
Our interpretation is that in maximization the boundary components are relatively faraway from each other and hence high accuracy results can be obtained for moderate values of $n$ and few number of iterations.
\end{remark}

\section{Conclusions}

Maximizing the conformal capacity of a constellation is opposite to minimizing studied in~\cite{hnv}. In~\cite{hnv} the main result was that the disks of the constellation group together in the local minima cases. Here we have shown that in the case of maximization the expected natural dispersion phenomenon occurs: the disks move as close to the unit circle as the constraints permit and, at the same time, the disks keep as far away from each other as possible. Replacing disks by other simple geometric objects also seems possible as our experiments with radial and rectilinear segments show.

A mathematical proof of the extremal cases we found in the experiments is missing. 
However, based on the numerical experiments we can predict the maximal capacities
of all constellations exhibiting sectorial symmetry.
In this case 
the constellation capacity is additive, equal to the sum of the capacities of the sectors or compartments. 
This was studied in~\cite{bsv} and~\cite{knv} from another point of view and similar conclusions obtained.

The study of this topic seems to offer many opportunities for later research. 
For example, one could study the above problems replacing the unit disk by some other domain, 
e.g. by a polygonal domain. Also one could investigate similar problems for other capacities 
such as the logarithmic and analytic capacities.

\section*{Acknowledgements} The authors would like to thank two anonymous reviewers for their valuable comments and suggestions which greatly improved the presentation of this paper.

\appendix
\section{Code Availability}
In the interest of reproducibility, 
the codes for our computations are available through the link \url{https://github.com/mmsnasser/maxcap}.

\bibliographystyle{elsarticle-num-names}

\end{document}